\documentclass[MSNbibl,number,seceqn,citesort,dvips]{arxbj}
\usepackage{upgreek}
\usepackage[left]{vmkcol}

\aid{0}
\volume{18}
\issue{4}
\pubyear{2012}
\firstpage{1421}
\lastpage{1447}
\doi{10.3150/11-BEJ384}

\makeatletter
\def\bsuffix #1{#1}

\newcolumntype{d}[1]{D{.}{.}{#1}}
\newcommand{\pn}{\stackrel{P}{\longrightarrow}}
\newcommand{\stab}{\stackrel{\mathcal{D}_{st}}{\longrightarrow}}
\newcommand{\Corr}{\operatorname{Corr}}

\newcommand{\cal}{\mathcal}
\newtheorem{theorem}{Theorem}
\newtheorem{lemma}{Lemma}
\newcommand{\eqref}[1]{(\ref{#1})}
\newcommand{\fraca}[2]{{#1}/{#2}}
\newcommand{\fracb}[2]{{(#1)}/{#2}}
\makeatother

\begin{document}
\begin{frontmatter}

\title{Model checks for the volatility under microstructure noise}
\runtitle{Model checks for the volatility}

\begin{aug}
\author{\fnms{Mathias} \snm{Vetter}\thanksref{e1}\ead[label=e1,mark]{mathias.vetter@rub.de}} \and
\author{\fnms{Holger} \snm{Dette}\corref{}\thanksref{e2}\ead[label=e2,mark]{holger.dette@rub.de}}
\runauthor{M. Vetter and H. Dette} 
\address{Ruhr-Universit\"at Bochum,
Fakult\"at f\"ur Mathematik,
44780 Bochum,
Germany.\\ \printead{e1,e2}}
\end{aug}

\received{\smonth{11} \syear{2009}}
\revised{\smonth{1} \syear{2011}}

%
\begin{abstract}
We consider the problem of testing the parametric form of the
volatility for high frequency data. It is demonstrated that in the
presence of
microstructure noise commonly used tests do not keep the preassigned
level and are inconsistent. The concept of preaveraging is used to
construct new tests, which do not suffer from these drawbacks. These
tests are based on a Kolmogorov--Smirnov or Cram\'{e}r--von-Mises
functional of an
integrated stochastic process,
for which weak convergence to a (conditional)
Gaussian process is established. The finite sample properties of a bootstrap
version of the test are illustrated by means of a simulation study.
\end{abstract}

%
\begin{keyword}
\kwd{goodness-of-fit test}
\kwd{heteroscedasticity}
\kwd{microstructure noise}
\kwd{parametric bootstrap}
\kwd{stable convergence}
\end{keyword}

\end{frontmatter}

\section{Introduction}\label{sec1}

The volatility is a popular measure of risk in finance with numerous
applications including the construction of optimal portfolios, hedging and
pricing of options. Therefore, estimating and investigating the
volatility and its dynamics is of particular importance in applications and
numerous models have been proposed for this purpose (see, e.g., Black
and Scholes~\cite{bs}, Vasicek~\cite{vas}, Cox \textit{\textit{et al.}}
\cite{cir}, Hull and White~\cite{hull} and Heston~\cite{heston} among
many others). Because the misspecification of the form of the volatility
can lead to serious consequences in the subsequent data analysis
numerous authors recommend to use goodness-of-fit tests for the postulated
model (see, e.g., Ait-Sahalia~\cite{ait}, Corradi and White \cite
{corr}, Dette \textit{et al.}~\cite{dpv}, Dette and Podolskij~\cite{dp} among
others).

In the present paper, we consider statistical inference in the case of
high frequency data, where for an increasing sample size information
about the whole path of the volatility is in principle available.
However, in concrete applications the situation is more complicated
because of the
presence of microstructure noise, which is usually persistent in such
data. This additional noise is caused by many sources of the
trading process such as discreteness of observations (see, e.g., Harris~\cite{harris1990},~\cite{harris1991}), bid-ask bounces or special
properties of the trading mechanism (see, e.g., Black~\cite{black1986}
or Amihud and Mendelson~\cite{am1987}). While microstructure noise has
been taken into account for the construction of estimators of the
integrated volatility and other related quantities (see, e.g., Zhang
\textit{et al.}~\cite{zha}, Jacod \textit{et al.}~\cite{jlmpv} or
Podolskij and Vetter~\cite{pv1},~\cite{pv2}), properties of goodness-of-fit
tests in this context have not been investigated so far in the literature.

\begin{table*}
\tabcolsep=0pt
\caption{Simulated level of the test (\protect\ref{1.2e}) for various
choices of $\omega$ and $\theta$, where the true
volatility function is $\sigma^2(t,x) = \theta+ (1-\theta)x^2$ and the
noise terms $U$ are normally distributed with mean zero and variance
$\omega^2$. In all cases, the sample size is given by $n= 16\,384$}
\label{tab1}\begin{tabular*}{\textwidth}{@{\extracolsep{\fill}}d{1.2}d{1.3}d{1.3}d{1.3}d{1.3}d{1.3}d{1.3}d{1.3}d{1.3}d{1.3}} \hline
\multicolumn{1}{@{}l}{$\omega$} & \multicolumn{3}{l}{0.01} &
\multicolumn{3}{l}{0.0025} & \multicolumn{3}{l@{}}{0.000625} \\[-5pt]
& \multicolumn{3}{l}{\hrulefill} & \multicolumn{3}{l}{\hrulefill} &
\multicolumn{3}{l@{}}{\hrulefill}\\
\multicolumn{1}{@{}l}{{$\theta$}/{$\alpha$}} & 0.025 & 0.05 & 0.1 &
0.025 & 0.05 & 0.1 & 0.025 & 0.05 & 0.1 \\
\hline
1 & 0.01 & 0.02 & 0.038 & 0.023 & 0.058 & 0.104 & 0.024 & 0.047 &
0.101\\
0.75 & 0.004 & 0.01 & 0.02 & 0.004 & 0.009 & 0.022 & 0.003 & 0.007 &
0.015\\
0.5 & 0.003 & 0.006 & 0.013 & 0.002 & 0.004 & 0.014 & 0.000 & 0.000 &
0.002\\
0.25 & 0.002 & 0.004 & 0.015 & 0.001 & 0.002 & 0.003 & 0.001 & 0.003 &
0.004\\
0 & 0.000 & 0.005 & 0.019 & 0.003 & 0.006 & 0.015 & 0.004 & 0.007 &
0.016\\
\hline
\end{tabular*}
\end{table*}

Consider for example the problem, where the process $\{Z_t\}_{t\in
[0,1]}$ is observed at the $n$ time points $0, 1/n, \ldots,1$. Under the
assumption that
$Z_t=X_t= \sigma_t \,\mathrm{d}W_{t}$,
Dette and Podolskij~\cite{dp} propose to reject the hypothesis of a
constant diffusion coefficient, that is, $H_0\dvt \sigma^2_t =
\sigma^2(t,X_{t}) = \sigma^2$, whenever
%
\begin{eqnarray}\label{1.2e}
\hspace*{-25pt}T_n (Z_1,\ldots, Z_n) &=& \sqrt{n} \sup_{t\in[0,1]}  \biggl| { \sum
_{k=1}^{\lfloor nt \rfloor} |Z_{\fraca kn} - Z_{\fracb{k-1}{n}}|^2 -
t \sum_{k=1}^{ n } |Z_{\fraca kn} - Z_{\fracb{k-1}{n}}|^2 \over
\sqrt{2} \sum_{k=1}^{ n } |Z_{\fraca kn} - Z_{\fracb{k-1}{n}}|^2 }
\biggr|  \nonumber
\\[-8pt]
\\[-8pt]
\hspace*{-25pt} &>& c_{1-\alpha} ,
\nonumber
\end{eqnarray}
where $c_{1-\alpha}$ denotes the $(1-\alpha)$-quantile of the supremum
of a Brownian Bridge. Now
consider the situation, where microstructure noise is present, which
is usually modeled by an additional additive component, that
is
%
\begin{equation}\label{1.3e}
Z_{\fraca in } =X_{\fraca in } + U_{\fraca in } ,  \qquad
i=1,\ldots,n,
\end{equation}
where $\{ U_{\fraca in } \mid  i=1,\ldots,n \}$ denotes a triangular array
of independent random variables
with mean $0$ and variance $\omega^2$.
In Table~\ref{tab1}, we show the finite sample behaviour of the test (\ref
{1.2e}) for the hypothesis of a constant volatility if
$\sigma^2_t = \sigma^2(t,x)=\theta+ (1- \theta)x^2$ (note that the
case $\theta=1$ corresponds to the null hypothesis).
We observe that the test keeps its preassigned level only in the case
where $\omega$ is rather small. In most cases, the nominal
level is clearly underestimated. On the other hand, the test is not
able to detect any alternative. An intuitive explanation
for this behaviour is that in the presence of microstructure noise the
increments
$Z_{\fraca in} - Z_{\fracb{i-1}{n}} = U_{\fraca in} - U_{\fracb{i-1}{n}} + \mathrm{O}_p(1/n)$
are dominated by the noise variables. This leads to inconsistent
estimates of the integrated volatility as pointed out in Zhang {\textit
{et al.}~\cite{zha}. More precisely, a straightforward calculation
shows that under microstructure noise the statistic $T_n (Z_1,\ldots,
Z_n) $ shows the same asymptotic behavior as
the statistic $ T_n (U_1,\ldots, U_n)$, which converges weakly\vadjust{\goodbreak}
to $\sqrt{\lambda/ 2} \sup_{t\in[0,1]} |B_t |, $
no matter if the null hypothesis is valid or not. Here $B_t$ denotes a
Brownian bridge and $\lambda= E[ (U_{ k/n} / \omega)^4 ]$. This means
that in the presence of microstructure noise the test
\eqref{1.2e} has asymptotic level $\alpha$ if and only if $\lambda=2$.
In all other cases, the test does not keep its preassigned level.
Moreover, because the asymptotic properties under null hypothesis and
alternative are the same, the test is not consistent.

The present paper is devoted to the problem of constructing a
consistent asymptotic level $\alpha$ test for a general parametric form
of the
volatility in the presence of microstructure noise. In Sections~\ref{sec2} and~\ref{sec3}, we present the basic model and introduce a stochastic process which
can be
used to test parametric hypotheses about the form of the volatility in
a noisy framework. Our main results are presented in Section~\ref{sec4}, where
we prove stable convergence of two such processes which form the basis
of the proposed goodness-of-fit tests. Section~\ref{sec5} deals with the problem
of testing nonlinear hypotheses for the volatility, whereas in Section~\ref{sec6} the finite sample properties of a bootstrap version of the new tests
are investigated. All proofs of the results are presented in the \hyperref[appm]{Appendix}.

\section{Testing parametric hypotheses for the volatility}\label{sec2}

Suppose that the process $X = (X_t)_t$ admits the representation
%
\begin{equation}\label{2.3}
X_{t}=X_0+ \int_0^t a_s \,\mathrm{d}s+ \int_0^t \sigma_s \,\mathrm{d}W_{s},
\end{equation}
where $W = (W_t)_t$ is a standard Brownian motion and the drift process
$a$ and the volatility process $\sigma$ satisfy some weak regularity
conditions, which will be specified later. Furthermore, we assume that
the process can be observed at discrete points on a fixed time interval,
say $[0,1]$.

Various assumptions on the structure of the volatility process have
been proposed in the literature. Among such models, a large class involves
the case where $\sigma$ is defined to be a local volatility process,
thus merely a function of time and state (see, e.g., Black and Scholes
\cite{bs}, Vasicek~\cite{vas}, Cox \textit{et al.}~\cite{cir}, Chan
\textit{et al.}~\cite{chan}, Ait-Sahalia~\cite{ait} or Ahn and Gao
\cite{ag} among many others). Because an appropriate modeling of the
volatility is of particular importance for the construction of portfolios,
hedging and pricing, many authors point out that the postulated model
should be validated by an appropriate goodness-of-fit test (see, e.g.,
Ait-Sahalia~\cite{ait} or Corradi and White~\cite{corr}). In several
cases, the hypothesis for the parametric form of the volatility is linear
and one has to consider the following two situations:
%
\begin{eqnarray} \label{h0} &&H_0\dvt \sigma_t^2 = \sigma^2(t,X_t) = \sum
_{i=1}^d \theta_i  \sigma_i^2(t,X_t) \qquad \forall t   \mbox{ a.s.} \quad\mbox{or}  \nonumber
\\[-9pt]
\\[-9pt]&&\bar
{H}_0\dvt \sigma_t = \sigma(t,X_t) = \sum_{i=1}^d \bar{\theta}_i  \bar{\sigma}_i(t,X_t) \qquad \forall t   \mbox{ a.s.},
\nonumber
\end{eqnarray}
where the functions $\sigma_1, \ldots, \sigma_d$ (or $\bar{\sigma}_1,
\ldots,
\bar{\sigma}_d$) are known and the parameters $\theta_1, \ldots, \theta
_d$ (or $\bar{\theta}_1, \ldots, \bar{\theta}_d$) are unknown, but
assumed to ensure\vadjust{\goodbreak} $\sigma^2(t,X_t) \geq0$ (or $\sigma(t,X_t) \geq0$)
almost surely. Other models involve volatility functions, where the
parameters enter nonlinearly (see Ait-Sahalia~\cite{ait}) and the
corresponding hypotheses will be considered later in Section~\ref{sec5}, because the
basic concepts are easier to explain in the linear context.

Let us focus on the problem involving $H_0$ for the moment, as the
other testing problem can be treated in the same way. Dette
and Podolskij~\cite{dp} propose to construct a test statistic using an
empirical version of the stochastic process
\begin{eqnarray*}
N_t = \int_0^t \Biggl \{ \sigma_s^2 - \sum_{j=1}^d \theta_j^{\mathrm{min}}  \sigma
_j^2 (s,X_s)  \Biggr\} \,\mathrm{d}s, \qquad\theta^{\mathrm{min}} = \operatorname{\arg\min}\limits_{\theta
\in\mathbb{R}^d} \int_0^1 \Biggl \{ \sigma_s^2 -
\sum_{j=1}^d \theta_j  \sigma_j^2 (s,X_s)  \Biggr\}^2 \,\mathrm{d}s.
\end{eqnarray*}
Thus, one uses the $L^2$ distance to determine the best approximation
to the unknown volatility process $\sigma^2$ by a linear combination of
the given functions $\sigma_1^2, \ldots, \sigma_d^2$. It can easily be
seen that $H_0$ is equivalent to $N_t = 0 \  \forall t    \mbox{ a.s.}$, and
a well-known result from Hilbert space theory (see Achieser~\cite{ach})
implies
%
\begin{eqnarray} \label{decomp} \theta^{\mathrm{min}} = D^{-1}C, \quad\mbox{thus} \quad N_t = B_t^0 - B_t^T D^{-1} C,
\end{eqnarray}
where
\begin{eqnarray*} B_t^0 = \int_0^t \sigma^2_s \,\mathrm{d}s \quad\mbox{and}\quad B_t^i = \int_0^t \sigma_i^2 (s, X_s) \,\mathrm{d}s  \qquad  \mbox{for }
 i=1,
\ldots, d,
\end{eqnarray*}
and $D$ and $C$ denote a $d \times d$-matrix and a $d$-dimensional
vector, respectively, with
\begin{eqnarray*} D_{ij} = \int^1_0 \sigma_i^2(s,X_s)  \sigma
_j^2(s,X_s) \,\mathrm{d}s \quad\mbox{and} \quad C_i = \int^1_0 \sigma_s^2  \sigma_i^2(s,X_s) \,\mathrm{d}s.
\end{eqnarray*}

In practice, one does not observe the entire path of the diffusion
process $X = (X_t)_t$ and it is therefore necessary to define an
empirical version based on appropriate estimators for the quantities in
(\ref{decomp}). Let us briefly discuss the solution to the problem in
the case, where $X$ can be observed without further
restrictions. Based on the decomposition above, Dette and Podolskij
\cite{dp} propose to define an empirical version $\tilde N_t = \tilde
B_t^0 - \tilde B_t^T \tilde D^{-1} \tilde C$, where one uses a
Riemann approximation of each integral, choosing $n |X_{\fraca kn} -
X_{\fracb{k-1}{n}}|^2$ as a local estimate for
$\sigma^2_{\fracb{k-1}{n}}$. Thus,
%
\begin{eqnarray} \label{dt}
\tilde{D}_{ij} &=& \frac1n \sum_{k=1}^n \sigma_i^2 \biggl({\frac kn},
X_{\fraca kn}\biggr)  \sigma_j ^2\biggl({\frac kn}, X_{\fraca kn}\biggr)  \qquad  \mbox
{for }  i,j=1, \ldots, d, \\
\tilde{C}_i & =&  \sum_{k=1}^n \sigma_i ^2\biggl({\frac
{k-1}{n}},X_{\fracb{k-1}{n}}\biggr)  |X_{\fraca kn} - X_{\fracb{k-1}{n}}|^2
 \qquad
\mbox{for }  i=1, \ldots, d,
\nonumber
\end{eqnarray}
and the quantities $\tilde B _t^0$ and $ \tilde{B}_t = ( \tilde
{B}^1_t, \ldots, \tilde{B}_t^d)^T$ are given
by
%
\begin{eqnarray} \label{2.21a}
  \tilde{B}_t^0 = \sum_{k=1}^{\lfloor nt \rfloor} |X_{\fraca kn} -
X_{\fracb{k-1}{n}}|^2, \qquad \tilde{B}_t^i = {\frac1n} \sum
_{k=1}^{\lfloor nt \rfloor} \sigma_i^2 \biggl({\frac kn}, X_{\fraca
kn}\biggr) \qquad  \mbox{for }  i=1, \ldots, d.
\end{eqnarray}
In this context, one can prove a (stable) central limit theorem for the
process $(\tilde N_t - N_t)_t$ with the optimal rate of convergence
$n^{-\fraca12}$, from which one may assess the distribution of suitable
test statistics. For example, if $d=1$, $\sigma_1^2(t,X_t)= 1$, the
hypothesis $H_0$ reduces to the hypothesis of constant volatility
considered in the introduction, and the Kolmogorov--Smirnov statistic
\eqref{1.2e} converges to the supremum of a Brownian bridge.

\section{Assumptions and definitions}\label{sec3} \label{assum}

Since we are dealing with microstructure noise, we have to define a
process $Z = (Z_t)_t$ which represents the noisy observations.
Typically one relates $Z$ to the underlying Ito
semimartingale $X$ through the equation $Z_t = X_t + U_t$ for some
noise process $U$. We restrict ourselves to the case of i.i.d. noise,
in which the process $U = (U_t)_t$ is independent of $X$ and satisfies
%
\begin{equation} \label{Umoment}
E[U_t]=0, \qquad E[U_t^2]=\omega^2, \qquad E[U_t^4] < \infty
\end{equation}
with a density having compact support. A precise definition of a proper
probability space that accommodates $Z$ can be found in Jacod \textit
{et al.}~\cite{jlmpv}. We assume further that $Z$ is observed at times
$0, 1/n, \ldots, 1$.
As pointed out in the introduction, the corresponding test based on
$\tilde N_t$ is not consistent for the hypothesis $H_0$ in the presence
of such microstructure noise. Thus, our aim is to define appropriate
estimators for the unknown quantities in (\ref{decomp}) in this noisy
framework, from which a more adequate statistic $\hat N_t$ can be
constructed. Note that in contrast to the previous setting we do not
only need a local estimator for the unknown volatility function $\sigma
^2$, but also for the (unobservable) path of $X$ itself.

The natural approach in order to construct estimators for the
volatility is to use increments of $Z$ as in the no-noise case, even though
a single increment does not provide sufficient information about $\sigma
^2$. This problem can be overcome by applying the idea of
pre-averaging, which was invented in Podolskij and Vetter~\cite{pv1}
and is based on moving averages of $Z$. To this end, we choose first a
sequence $m_n$, such that
%
\begin{eqnarray}\label{kappa} \frac{m_n}{\sqrt n} = \kappa+ \mathrm{o}(n^{-
\fraca14})
\end{eqnarray}
for some $\kappa> 0$, and a nonzero real-valued function $g\dvtx \mathbb
{R} \rightarrow
\mathbb{R}$, which vanishes outside of the interval $(0,1)$, is
continuous and piecewise $C^1$ and has a piecewise Lipschitz derivative~$g'$. We associate
with $g$ (and $n$) the following real valued numbers and functions:
%
\begin{equation} \label{psi}
\hspace*{-15pt}\everymath{\displaystyle }\cases{ g_j^n =
g\biggl(\frac{j}{m_n}\biggr),\qquad
g_j^{\prime n} = g_j^n - g_{j+1}^n, \qquad   \psi_1 = \int_0^1 (g'(s))^2
\,\mathrm{d}s,\qquad
 \psi_2 = \int_0^1 (g(s))^2 \,\mathrm{d}s, \cr s \in[0,1]  \mapsto \phi_1(s) = \int_s^1 g'(u) g'(u-s) \,\mathrm{d}u,  \qquad  \phi_2(s) = \int_s^1 g(u)
g(u-s) \,\mathrm{d}u, \cr
i,j=1,2{:} \qquad    \Phi_{ij} = \int_0^1 \phi_i(s) \phi_j(s)
\,\mathrm{d}s.
}
\end{equation}
Finally, we define for an arbitrary process $V$ the preaveraged statistic
%
\begin{eqnarray}
\overline{V}^n_k=\sum_{j=1}^{m_n} g^n_j \Delta_{k+j}^n V,
\end{eqnarray}
where $\Delta_j^n V = V_{\fraca jn} - V_{\fracb{j-1}{n}}$. Due to the
assumptions on $g$ the pre-averaged statistic $\overline{Z}^n_k$
reduces the impact of the noise, but still provides information about
the increments of $X$ (and thus locally about $\sigma$). Precisely, we have
%
\begin{eqnarray} \label{orders}
\overline{X}^n_k = \mathrm{O}_p\Biggl(\sqrt{\frac{m_n}{n}}\Biggr) \quad\mbox{and} \quad
\overline{U}^n_k = \mathrm{O}_p\Biggl(\sqrt{\frac{1}{m_n}}\Biggr),
\end{eqnarray}
and by definition of $m_n$ both terms are of the same order. This means
in particular that statistics based on $\overline{Z}^n_k$ are in
general biased when used for volatility estimation, but it turns out
that a larger choice of $m_n$ results in a worse rate of convergence.
See Podolskij and Vetter~\cite{pv1} for details.

An estimator for $X_{\fraca{k}{n}}$ can be constructed in a similar way:
We set
%
\begin{eqnarray} \label{hat} \hat{X}_{\fraca kn} = \frac{1}{m_n} \sum
_{j=1}^{m_n} Z_{\fracb{k+j}{n}},
\end{eqnarray}
and it is easy to see that this procedure reduces the impact of the
noise variables around time $\frac kn$, but still provides information
about the latent price
$X_{\fraca kn}$, since the path of $X$ is H\"older continuous of any
order $\alpha< 1/2$. Also one observes essentially from (\ref{orders})
that the auxiliary sequence $m_n$ is chosen in the optimal way, giving
the smallest possible size for the approximation error.

As pointed out before, we need additional assumptions on the process
$X$ as well as on the given basis functions in $H_0$ and $\bar H_0$,
respectively. Since the conditions on $\sigma_i^2$ and $\bar\sigma_i$
are similar, we will restrict ourselves to the first case
only.

It is required that the functions $\sigma^2_1, \ldots, \sigma^2_d$ are
linearly independent and that each $\sigma_i^2$ is twice continuously
differentiable. Moreover, we assume that $E[|\det(D)|^{-\beta}] < \infty
$ for some $\beta> 0$.

Regarding the various processes in $X$, the assumptions are as weak as
possible when testing for $H_0$. We simply have to ensure that the
process in (\ref{2.3}) is well defined, which follows if we assume that
$a$ is locally bounded and predictable and that $\sigma$ is c\`adl\`ag
(see Jacod and Shiryaev~\cite{jacshir} or Revuz and Yor~\cite{rev}).
When working with $\bar H_0$ we propose additionally that the true
volatility process $\sigma$ is almost surely positive and that is has a
representation of the form (\ref{2.3}) as well, namely that it
satisfies
\begin{eqnarray*} \sigma_t=\sigma_0+\int_0^t a'_s\,\mathrm{d}s+\int_0^t\sigma
'_{s}\,\mathrm{d}W_s +\int_0^t v'_{s}\,\mathrm{d}V_s,
\end{eqnarray*}
where $a'$, $\sigma'$ and $v'$ are
adapted c\`adl\`ag processes, with $a'$ also being predictable and
locally bounded, and $V$ is a second Brownian motion, independent of
$W$. Moreover, $a$ is supposed to be c\`agl\`ag.

\section{Goodness-of-fit tests addressing microstructure noise}\label{sec4}

We start with the construction of a test for the hypothesis $H_0$
again. Local estimators for the volatility can now be obtained from
$|\overline{Z}^n_k|^2$, but we have seen before that this quantity is
not an unbiased estimate for $\sigma_{\fraca kn}^2$ and that it has a
different stochastic order than the increments $X_{\fraca kn} -
X_{\fracb
{k-1}{n}}$ in the no-noise case. A corrected statistic (see Jacod
\textit{et al.}~\cite{jlmpv}) is given by
%
\begin{eqnarray} \label{sigmahat}
\hat\sigma^2_{\fraca kn} = \frac{n^{\fraca12}}{\kappa\psi_2}
\biggl(|\overline{Z}^n_k|^2 - n^{-\fraca12} \frac{\psi_1}{\kappa} \hat{\omega
}^2_n  \biggr)  \qquad \mbox{with }  \hat{\omega}^2_n = \frac{1}{2n} \sum
_{i=1}^n |\Delta_i^n Z|^2,
\end{eqnarray}
where the latter term is a consistent estimator for $\omega^2$, see Zhang
\textit{et al.}~\cite{zha}. Mimicking the procedure from the no-noise
case presented in Section~\ref{sec2}, we set
%
\begin{eqnarray} \label{dhat}
\hspace*{-25pt}\hat{D}_{ij} = \frac1n \sum_{k=1}^{n-m_n} \sigma_i^2 \biggl({\frac kn}, \hat
{X}_{\fraca kn}\biggr)  \sigma_j ^2\biggl({\frac kn}, \hat{X}_{\fraca kn}\biggr) \quad
\mbox{and} \quad\hat{C}_i = \frac{1}{n} \sum_{k=1}^{n-m_n} \sigma
_i^2\biggl({\frac kn},\hat{X}_{\fraca kn}\biggr)  \hat\sigma^2_{\fraca kn}
\end{eqnarray}
as well as
%
\begin{equation}\label{Bt0}
\hat{B}_t^0 = \frac{1}{n} \sum_{k=1}^{\lfloor nt \rfloor- m_n} \hat
\sigma^2_{\fraca kn} \quad\mbox{and} \quad\hat{B}_t^i = {\frac1n}
\sum_{k=1}^{\lfloor nt \rfloor-m_n} \sigma_i^2 \biggl({\frac kn}, \hat
{X}_{\fraca kn}\biggr)
\end{equation}
for $i,j=1, \ldots, d$. We define at last the process
%
\begin{eqnarray} \label{Nt} \hat N_t = \hat B_t^{0} - \hat B_t^T \hat
D^{-1} \hat C,
\end{eqnarray}
which turns out to be an appropriate estimate of the process $\{ N_t\}
_{t \in[0,1]}$. Our first result specifies the
asymptotic properties of the process $\{ A_n(t) \}_{t\in[0,1]}$ with
\mbox{$A_n(t) = n^{\fraca14} (\hat N_t - N_t).$}

\begin{theorem} \label{tigh}
If the assumptions stated in the previous sections are satisfied, the
process $(A_n(t))_{t \in[0,1]}$ converges weakly in $D[0,1]$ to a mean
zero process $(A(t))_{t \in[0,1]}$. Conditionally on $\mathcal F$ the
limiting process is Gaussian, and its finite dimensional distributions
coincide with the conditional (with respect to $\mathcal F$) finite
dimensional distributions of the process
%
\begin{equation} \label{covkernel}
\hspace*{-15pt} \biggl\{\! \gamma_V \bigl ( I \{ V \leq t \} - B_t^T D^{-1} h(V,X_V)  \bigr)
-  \biggl( \int_0^t\! \gamma_s \,\mathrm{d}s - B_t^T D^{-1}\! \int_0^1\! \gamma_s  h(s,X_s)
\,\mathrm{d}s  \biggr) \! \biggr\}_{t \in[0,1]}\!,
\end{equation}
where $V \sim{\cal U} [0,1]$, $h(s,X_s) =( \sigma_1^2 (s, X_s), \ldots
, \sigma_d^2 (s, X_s) )^T$ and
%
\begin{eqnarray} \label{gamma} \gamma_s^2 = \frac{4}{\psi_2^2} \biggl (\Phi
_{22} \kappa\sigma_s^4 + 2 \Phi_{12} \frac{\sigma_s^2 \omega^2}{\kappa
} + \Phi_{11}
\frac{\omega^4}{\kappa^3} \biggr).
\end{eqnarray}
%
\end{theorem}
We see from Theorem~\ref{klord} in the \hyperref[appm]{Appendix} that the asymptotics is
only driven by $\hat B_t^0$ and $\hat C$. The error due to the
estimation of $B_t$ and $D$ is of small order, which explains the
particular form of the limiting distribution. Note also that the rate of
convergence $n^{-\fraca14}$ is optimal for this problem, since it is
already optimal for the estimation of $B_t^0$ even in a parametric setting
(cf. Gloter and Jacod~\cite{glot}).

In order to construct a test statistic based on Theorem~\ref{tigh}, we
have to define an appropriate estimator for the conditional variance of
the process $\{ A(t)\}_{t \in[0,1]}$, which is given by
\begin{eqnarray*} s_t^2 = \int_0^t \gamma_s^2 \,\mathrm{d}s - 2 B_t^T D^{-1} \int
_0^t \gamma_s^2 g (s, X_s) \,\mathrm{d}s + B_t^T D^{-1} \int_0^1 \gamma_s^2 g (s,
X_s) g ^T(s, X_s) \,\mathrm{d}s   D^{-1}B_t.
\end{eqnarray*}
Obviously, we use $\hat{B}_t$ and $\hat D$ as the empirical
counterparts for $B_t$ and $D$. In order to obtain estimates for the other
random elements of $s_t^2$, note that $\gamma_s^2$ plays a key role in
Jacod \textit{et al.}~\cite{jlmpv} as well, where it is the (local)
conditional variance in a central limit theorem for $n^{1/4} (\hat
B_t^{0} - B_t^{0})$. Thus, in accordance to that paper we define
\begin{eqnarray*}  \Gamma_k &=& \frac{4  \Phi_{22}}{3  \kappa \psi
_2^4}  |\overline Z_k^n|^4 + n^{-\fraca12}  \frac{8}{\kappa^2}
 \biggl (\frac{\Phi_{12}}{\psi_2^3} -\frac{\Phi_{22}  \psi_1}{\psi_2^4} \biggr)  |\overline Z_k^n|^2  \hat{\omega}^2\\
 &&{} +
n^{-1} \frac{4}{\kappa^3} \biggl (\frac{\Phi_{11}}{\psi_2^2} - \frac{ 2  \Phi_{12}  \psi_1}{\psi_2^3}
 + \frac{\Phi_{22}  \psi_1^2}{\psi_2^4} \biggr)  \hat{\omega}^4,
\end{eqnarray*}
which is a local estimator for the process $\gamma^2$ after rescaling.
Thus, we set
\begin{eqnarray*}
\hat g_0(t) &=& \sum_{k=1}^{\lfloor nt \rfloor- m_n} \Gamma_k \pn\int
_0^t \gamma_s^2 \,\mathrm{d}s, \\
 g_i (t) &=& \sum_{k=1}^{\lfloor nt \rfloor-
m_n} \Gamma_k  \sigma_i^2 \biggl(\frac{k-1}{n}, \hat X_{\fracb{k-1}{n}}\biggr) \pn
\int_0^t \gamma_s^2  \sigma_i^2(s,X_s) \,\mathrm{d}s, \\
\hat g_{ij} &=& \sum_{k=1}^{n} \Gamma_k  \sigma_i^2 \biggl(\frac{k-1}{n},
\hat X_{\fracb{k-1}{n}}\biggr)  \sigma_j^2 \biggl(\frac{k-1}{n}, \hat
X_{\fracb
{k-1}{n}}\biggr) \pn\int_0^1 \gamma_s^2  \sigma_i^2(s,X_s)  \sigma
_j^2(s,X_s) \,\mathrm{d}s.
\end{eqnarray*}
Inserting these estimators into the corresponding elements of $s_t^2 $
gives the consistent estimator
%
\begin{eqnarray} \label{est} \hat s^2_t
= \hat g_0(t) - 2 \hat B_t^T \hat D^{-1} \hat g (t) + \hat B_t^T \hat
D^{-1} \hat G \hat D^{-1} \hat B_t,
\end{eqnarray}
where $\hat g (t) = (\hat
g_1 (t), \ldots, \hat g_d (t))^T$ and $\hat G =( \hat
g_{ij})_{i,j=1}^d$. A consistent test for the hypothesis $H_0$ is now
obtained by
rejecting the null hypothesis for large values of Kolmogorov--Smirnov or
Cram\'{e}r--van-Mises functional of the process
$
\{{n^{1/4} \hat N_t} /{ {\hat s_t}} \}_{t \in[0,1]}.
$
Note however that the distribution of this process is not feasible in
general: even though for each fixed $t$ the statistic $n^{1/4} \hat N_t
/ {\hat s_t}$ converges weakly to a standard normal distribution, the
covariance structure of the process typically depends on the entire
(unobservable) process $(X_t)_t$. For this reason, we will later use a
bootstrap procedure to obtain critical values.\vadjust{\goodbreak}

In principle, a similar approach can be used to construct a test for
the hypothesis $\bar H_0$. However, in this case things change
considerably. Dette and Podolskij~\cite{dp} restate this hypothesis as
$M_t = 0$    $\forall t    \mbox{ a.s.}$, where
%
\begin{eqnarray}\label{2.17b}
M_t &=& \int_0^t  \Biggl\{ \sigma_s - \sum_{j=1}^d \bar{\theta}_j^{\mathrm{min}} \bar
{\sigma}_j (s,X_s)  \Biggr\} \,\mathrm{d}s,
 \nonumber
\\[-9pt]
\\[-9pt]
\bar{\theta}^{\mathrm{min}}&=&
\operatorname{\arg\min}\limits_{ \bar{\theta} \in\mathbb{R}^d} \int_0^1  \Biggl\{ \sigma_s -
\sum_{j=1}^d \bar{\theta}_j \bar{\sigma}_j (s,X_s)  \Biggr\}^2 \,\mathrm{d}s.
\nonumber
\end{eqnarray}
Obviously, we have an analogous representation as in (\ref{decomp}),
namely \(M_t = R_t^0 - R_t^T Q^{-1} S\), where
\begin{eqnarray*} R_t^0 = \int_0^t \sigma_s \,\mathrm{d}s \quad\mbox{and}
\quad R_t^i = \int_0^t \bar
\sigma_i (s, X_s) \,\mathrm{d}s  \qquad  \mbox{for }  i=1, \ldots, d,
\end{eqnarray*}
and $Q$ and $S$ are a $d \times d$-matrix and a $d$-dimensional vector,
respectively, with
\begin{eqnarray*} Q_{ij} = \int^1_0 \bar\sigma_i(s,X_s)  \bar\sigma
_j(s,X_s) \,\mathrm{d}s \quad\mbox{and} \quad S_i = \int^1_0 \sigma_s  \bar
\sigma_i(s,X_s) \,\mathrm{d}s.
\end{eqnarray*}
However, an appropriate definition of an empirical version of the form
$\hat M_t = \hat R_t^0 - \hat R_t^T \hat
Q^{-1} \hat S$ requires some less obvious modifications, because
local estimators for $\sigma_s$ are more difficult to obtain in this
setting. Using a preaveraged estimator of the form
$|\overline{Z}^n_k|$ again causes an intrinsic bias, but due to the
absolute value (instead of the square as in the previous setting) its
correction turns out to be impossible at the optimal rate. However, we
can see from (\ref{orders}) that using in
(\ref{kappa}) a sequence of a larger magnitude than $n^{\fraca12}$
reduces the impact of the noise terms in $\overline Z_k^n$. This
modification makes inference about $\sigma_s$ possible, though
resulting in a worse rate of convergence. To be precise, we fix some
$\delta>
\frac16$ and choose $l_n$ such that
\begin{eqnarray*} \frac{l_n}{n^{\fraca12 + \delta}} = \rho+
\mathrm{o}\bigl(n^{-(\fraca14 + \fraca{\delta}{2})}\bigr)
\end{eqnarray*}
for some $\rho> 0$. Using the sequence $l_n$ instead of $m_n$, we
define all quantities from (\ref{psi}) to (\ref{hat}) in the
straightforward way. Next, we set
\begin{eqnarray*}
\bar\sigma_{\fraca kn} = n^{\fraca14 - \fraca{\delta}{2}} \frac{1}{\sqrt
{\rho\psi_2} \mu_1} |\overline{Z}^n_k|
\end{eqnarray*}
as a local estimator for $\sigma_{\fraca kn}$, where $\mu_1$ denotes the
first absolute moment of a standard normal distribution. In a similar
way as before,
\begin{eqnarray*}
\hat{Q}_{ij} = \frac1n \sum_{k=1}^{n-l_n} \bar\sigma_i \biggl({\frac kn},
\hat{X}_{\fraca kn}\biggr)  \bar\sigma_j\biggl({\frac kn}, \hat{X}_{\fraca kn}\biggr)
\quad\mbox{and} \quad\hat{S}_i = \frac{1}{n} \sum_{k=1}^{n-l_n} \bar
\sigma_i\biggl({\frac kn},\hat{X}_{\fraca
kn}\biggr)  \bar\sigma_{\fraca kn}\vadjust{\goodbreak}
\end{eqnarray*}
as well as
\begin{eqnarray*}
\hat{R}_t^{0} = \frac{1}{n} \sum_{k=1}^{\lfloor nt \rfloor- l_n} \bar
\sigma_{\fraca kn} \quad\mbox{and} \quad\hat{R}_t^i = {\frac1n} \sum
_{k=1}^{\lfloor nt
\rfloor-l_n} \bar\sigma_i \biggl({\frac kn}, \hat{X}_{\fraca kn}\biggr)
\end{eqnarray*}
for $i,j=1, \ldots, d.$ Finally, we define $B_n(t) = n^{\fraca14 -
\fraca\delta2} (\hat M_t - M_t)$ for any $t \in[0,1]$ and obtain the
following result.

\begin{theorem} \label{absol}
If the assumptions stated in the previous sections are satisfied, the
process $(B_n(t))_{t \in[0,1]}$ converges weakly in $D[0,1]$ to a mean
zero process $(B(t))_{t \in[0,1]}$. Conditionally on $\mathcal F$ the
limiting process is Gaussian, and its finite dimensional distributions
coincide with the conditional (with respect to $\mathcal F$) finite
dimensional distributions of the process
%
\begin{equation} \label{covkerne}
\hspace*{-15pt} \biggl\{\! \bar\gamma_V  \bigl( I \{ V \leq t \} - R_t^T Q^{-1} \bar
h(V,X_V)  \bigr) -  \biggl( \int_0^t\! \bar\gamma_s \,\mathrm{d}s - R_t^T Q^{-1}\! \int
_0^1\! \bar\gamma_s   \bar h(s,X_s) \,\mathrm{d}s  \biggr)  \!\biggr\}_{t \in[0,1]}\!,
\end{equation}
where $V \sim{\cal U} [0,1]$, $ \bar h(s,X_s) =( \bar\sigma_1 (s,
X_s), \ldots, \bar\sigma_d (s, X_s) )^T $ and
%
\begin{eqnarray} \label{f}
\bar\gamma_s^2 &=& \frac{2 \rho\Xi}{\mu_1^2} \sigma_s^2 , \qquad\Xi=
\int_0^1 \xi(s) \,\mathrm{d}s,   \qquad   \xi(s) = f \biggl(\frac{\phi_2(s)}{\psi_2} \biggr),
 \nonumber
 \\[-8pt]
 \\[-8pt] f(u) &=& \frac2 \uppi \bigl(u \arcsin(u) + \sqrt{1- u^2} -1 \bigr).
\nonumber
\end{eqnarray}
\end{theorem}

The estimation of the conditional variance of the process $\{ B(t)\}_{t
\in[0,1]}$,
\begin{eqnarray*} r_t^2 = \int_0^t \bar \gamma_s^2 \,\mathrm{d}s
- 2 R_t^T Q^{-1} \int_0^t \bar \gamma_s^2 \bar h (s, X_s) \,\mathrm{d}s +
R_t^T Q^{-1} \int_0^1 \bar\gamma_s^2 \bar h (s, X_s) \bar g ^T(s, X_s)
\,\mathrm{d}s   Q^{-1}R_t,
\end{eqnarray*}
becomes easier in this context, as the order of $l_n$ is chosen in such
a way that no characteristics of $U$ are involved anymore. A natural
estimator for $\sigma_{\fraca kn}^2$ becomes
\begin{eqnarray*} \bar\Gamma_k = n^{-(\fraca12 + \delta)}  \frac{2  \Xi}{\psi_2  \mu_1^2}  |\overline Z_k^n|^2,
\end{eqnarray*}
thus
\begin{eqnarray*}
\hat h_0 (t) &=& \sum_{k=1}^{\lfloor nt \rfloor- l_n} \bar\Gamma_k
\pn\int_0^t \bar{\gamma}_s^2 \,\mathrm{d}s,  \\
\hat h_i (t) &=& \sum
_{k=1}^{\lfloor nt \rfloor- l_n} \bar\Gamma_k  \bar\sigma_i \biggl(\frac
{k-1}{n}, \hat X_{\fracb{k-1}{n}}\biggr) \pn\int_0^t \bar\gamma_s^2  \bar
\sigma_i(s,X_s) \,\mathrm{d}s, \\
\hat h_{ij} &=& \sum_{k=1}^{n} \bar\Gamma_k  \bar\sigma_i \biggl(\frac
{k-1}{n}, \hat X_{\fracb{k-1}{n}}\biggr)  \bar\sigma_j \biggl(\frac{k-1}{n}, \hat
X_{\fracb{k-1}{n}}\biggr) \pn\int_0^1 \bar\gamma_s^2  \bar\sigma_i(s,X_s)
 \bar\sigma_j(s,X_s) \,\mathrm{d}s,
\end{eqnarray*}
and consequently a consistent estimator
$\hat r^2_t$ for the conditional variance is given by
%
\begin{eqnarray} \label{estabs} \hat r^2_t = \hat h_0(t) - 2 \hat R
_t^T \hat Q^{-1} \hat h (t) +
\hat R_t^T \hat Q^{-1} \hat H \hat Q^{-1} \hat R_t,
\end{eqnarray}
where $\hat h (t) = (\hat h_1 (t), \ldots, \hat h_d (t))^T$ and $\hat
H =( \hat h_{ij})_{i,j=1}^d$.
A consistent test for the hypothesis $\bar H_0$ is now obtained by
rejecting the null hypothesis for large values of the Kolmogorov--Smirnov
or Cram\'{e}r--van-Mises functional of the process
$
\{ {n^{1/4 - \delta/2} \hat M_t}/ { {\hat r_t}} \}_{t \in[0,1]}.
$

Note that one knows from previous work that it is neither necessary to
define $X$ to be an Ito semimartingale with continuous paths as in
\eqref{2.3} nor to model the noise terms $U$ as being independent and
identically distributed to obtain similar results as in Theorems~\ref{tigh} and~\ref{absol}.
In fact, for an underlying Ito
semimartingale exhibiting jumps one can use bipower-type estimators as
discussed in
Podolskij and Vetter~\cite{pv2} in order to define an estimator closely
related to $\hat B_t^0$. Moreover, it has been argued in Jacod
\textit{et al.}~\cite{jlmpv} that even for a noise process with a c\`
adl\`ag variance a similar theory as presented in this paper applies.

\section{Nonlinear hypotheses}\label{sec5}

In this section, we briefly discuss the case of a nonlinear hypothesis
%
\begin{eqnarray} \label{5.1} H_0\dvt \sigma^2_t = \sigma^2(t,X_t)=\sigma
^2(t, X_t,
\theta) \qquad \forall t   \mbox{ a.s.},
\end{eqnarray}
where $\theta\in\Theta\subset\mathbb{R}^d$ denotes the unknown
parameter and $\sigma^2$ satisfies some differentiability assumption.
As before, we restate $H_0$ as $N_t = 0 \  \forall t    \mbox{ a.s.}$, where
$N_t$ is the difference between the true integrated volatility and its
best $L^2$-approximation from the parametric class. Therefore, we set $
N_t = B_t^0 - B_t (\theta_0)$ with $B_t^0$ from above and $B_t (\theta)
= \int_0^t \sigma^2(s, X_s, \theta) \,\mathrm{d}s$. We have $\theta_0 = \operatorname{\arg\min}_{\theta\in\Theta}  f(\theta)$ with
\begin{eqnarray*}
f(\theta) = \int_0^t  \{ \sigma_s^2 - \sigma^2(s, X_s, \theta)
\}^2 \,\mathrm{d}s.
\end{eqnarray*}
In order to obtain some $\hat N_t$, we use $\hat B^0_t$ from (\ref
{Bt0}) and need estimates for $B_t (\theta)$ and $f(\theta)$. We set
%
\begin{eqnarray} \label{5.3}
\hat B_t (\theta) &=& \frac{1}{n} \sum^{\lfloor nt \rfloor-m_n}_{k=1}
\sigma^2 \biggl(\frac{k}{n}, \hat X_{\fraca{k}{n}},\theta\biggr) \quad\mbox{and}\nonumber\\[-8pt]\\[-8pt]
 f_n(\theta) &=& \frac1n \sum^{n-m_n}_{k=1}  \biggl\{ \hat\sigma
^2_{\fraca kn} - \sigma^2 \biggl(\frac{k}{n}, \hat X_{\fraca{k}{n}}, \theta\biggr)
 \biggr\}^2,\nonumber
\end{eqnarray}
and with $\hat\theta= \operatorname{\arg\min}_{\theta\in\Theta}  f_n(\theta
)$ we define $\hat N_t = \hat B_t^0 - \hat B_t (\hat\theta)$.

When deriving the asymptotic distribution of $n^{\fraca14} (\hat N_t -
N_t)$, the difference compared to the previous section regards only
$\hat B_t (\theta_0) - B_t (\hat\theta)$. In the following, we will
give some hints that explain why that discrepancy is actually quite
small. In fact, we will show that
%
\begin{eqnarray} \label{5.6}
\hspace*{-15pt}\hat B_t (\hat\theta) - B_t (\theta_0) = - \int_0^t  \biggl( \frac
{\partial}{\partial\theta} \sigma^2 (s, X_s, \theta) \bigg |_{\theta
=\theta_0} \biggr)^{T}\,\mathrm{d}s  (f''(\theta_0))^{-1}  f_n' (\theta_0) +
\mathrm{o}_p(n^{-\fraca14})
\end{eqnarray}
holds. Thus there is a one-to-one correspondence to the linear case, as
the first two quantities are analogues of $B_t^T$ and $D^{-1}$, whereas
$-f_n' (\theta_0)$ plays the role of $\hat C - C$. Consequently, the
process $n^{\fraca14}(\hat N_t - N_t)$ exhibits a similar asymptotic
behavior as in the linear case.

In order to prove (\ref{5.6}), note from similar arguments as in the
proof of Theorem~\ref{klord} that
%
\begin{eqnarray} \label{5.2}
\hat B_t (\hat\theta) - B_t (\theta_0) = \int_0^t  \{ \sigma^2(t,
X_t, \hat\theta) - \sigma^2(t, X_t, \theta_0)  \} \,\mathrm{d}s +
\mathrm{o}_p(n^{-\fraca14}).
\end{eqnarray}
Under common regularity conditions for nonlinear regression (see
Gallant~\cite{gallant1987} or Seber and Wild~\cite{sw1989}), $\theta_0$
is the unique minimum of $f$ and attained at an interior point of
$\Theta$. It is easy to see that $\hat\theta\rightarrow\theta_0$ in
probability in this case, and thus we can assume that $\hat\theta$
satisfies $f_n' (\hat\theta) = 0$. This implies
\begin{eqnarray*}
0 = f_n' (\hat\theta) = f_n' (\theta_0) + f''_n(\tilde\theta)(\hat
\theta- \theta_0)   \quad \Leftrightarrow  \quad \hat\theta- \theta_0 =
- (f''_n(\tilde\theta))^{-1}  f_n' (\theta_0)
\end{eqnarray*}
for an appropriate choice of $\tilde\theta$. We have $\tilde\theta
\rightarrow\theta_0$ in probability as well, and therefore it can be
assumed that the $d \times d$-dimensional matrix $f''_n(\tilde\theta)$
is positive definite and that the difference $\|f''_n(\tilde\theta) -
f''_n(\theta_0)\|$ is small. Furthermore, $f''_n (\theta_0)$ takes the form
\begin{eqnarray*}
f''_n (\theta_0) &=& 2 \Biggl ( \frac1n S^T S - \frac1n \sum
^{n-m_n}_{k=1}  \biggl\{ \hat\sigma^2_{\fraca kn} - \sigma^2 \biggl(\frac
{k}{n}, \hat X_{\fraca{k}{n}}, \theta_0\biggr)  \biggr\} H_k  \Biggr),
\end{eqnarray*}
where the $(n-m_n) \times d$ matrix $S$ and the Hessian $H_k$ are given by
\begin{eqnarray*}
S=  \biggl( \frac{\partial}{\partial\theta} \sigma^2 \biggl(\frac{k}{n},
\hat X_{\fraca{k}{n}}, \theta\biggr)  \bigg |_{\theta=\theta_0}
\biggr)_{k=1,\ldots,n-m_n} \quad\mbox{and} \quad H_k = \frac{\partial
^2}{\partial\theta^2} \sigma^2\biggl(\frac kn,X_{\fraca kn},\theta\biggr)
\bigg|_{\theta=\theta_0}.
\end{eqnarray*}
From the same arguments that lead to \eqref{5.2}, we have $f''_n (\theta
_0) = f'' (\theta_0) + \mathrm{O}_p(n^{-\fraca14})$, where
\begin{eqnarray*}
f'' (\theta_0) &=& 2 \int_0^1  \biggl( \frac{\partial}{\partial\theta}
\sigma^2 (s, X_s, \theta)   \bigg|_{\theta=\theta_0}  \biggr)^T  \biggl(
\frac{\partial}{\partial\theta} \sigma^2 (s, X_s, \theta)
\bigg|_{\theta=\theta_0}  \biggr)\,\mathrm{d}s \\
&&{}- 2 \int_0^1  \biggl\{ \sigma_s^2 -
\sigma^2(s, X_s, \theta_0)  \biggr\} \frac{\partial^2}{\partial\theta^2}
\sigma^2(s,X_s,\theta) \bigg|_{\theta=\theta_0}\,\mathrm{d}s
\end{eqnarray*}
is positive definite. Note that the second term in this sum vanishes,
when either the hypothesis is linear (since the Hessian is zero) or the
null hypothesis is valid (since $\sigma^2_s$ equals $\sigma^2(s,
X_s,\theta_0)$). In these cases the matrix $f'' (\theta_0)$ takes
precisely the same form as $D$ in the linear setting. In any case, $
f'' (\theta_0)$ is of order $\mathrm{O}_p(1)$.

Regarding $f_n' (\theta_0)$, a similar calculation as given in the
\hyperref[appm]{Appendix} plus the definition of $\theta_0$ yield
\begin{eqnarray*}
-f_n' (\theta_0) &=&
2  \Biggl( \frac1n \sum^{n-m_n}_{k=1} \hat\sigma^2_{\fraca kn}  \frac
{\partial}{\partial\theta} \sigma^2 \biggl(\frac{k}{n}, \hat X_{\fraca
{k}{n}}, \theta\biggr)  \bigg |_{\theta=\theta_0} - \int_0^1 \sigma_s^2  \frac{\partial}{\partial\theta} \sigma^2 (s, X_s, \theta)
\bigg|_{\theta=\theta_0} \,\mathrm{d}s  \Biggr)\\
&&{} + \mathrm{o}_p(n^{-\fraca14}),
\end{eqnarray*}
and thus $f_n' (\theta_0)$ is of order $\mathrm{O}_p(n^{-\fraca14})$, just as
$\hat C - C$. We conclude that $\hat\theta- \theta_0 =
\mathrm{O}_p(n^{-\fraca
14})$ as well, and a Taylor expansion gives (\ref{5.6}).

\begin{table*}[b]
\tabcolsep=0pt
\tablewidth=310pt
\caption{Simulated level of the bootstrap test proposed by Dette and
Podolskij \protect\cite{dp}, where the volatility function equals $
H_0\dvt
\sigma^2(t,x) = \theta x^2$, but the observations are corrupted with
normally distributed noise having variance $\omega^2$}
\label{tab2}\begin{tabular*}{310pt}{@{\extracolsep{\fill
}}d{1.3}d{1.3}d{1.3}d{1.3}d{1.3}d{1.3}d{1.3}@{}} \hline
\multicolumn{1}{@{}l}{$n$} & \multicolumn{3}{l}{256} & \multicolumn
{3}{l@{}}{1024}
\\[-5pt]
& \multicolumn{3}{l}{\hrulefill} & \multicolumn{3}{l@{}}{\hrulefill}\\
\multicolumn{1}{@{}l}{{$\omega$}/{$\alpha$}} & 0.025 & 0.05 & 0.1 &
0.025 & 0.05 & 0.1 \\
\hline
0.001 & 0.033 & 0.062 & 0.111 & 0.333 & 0.415 & 0.512\\
0.002 & 0.158 & 0.243 & 0.324 & 0.810 & 0.862 & 0.907 \\
0.004 & 0.392 & 0.518 & 0.650 & 0.993 & 0.996 & 0.998 \\
0.005& 0.497 & 0.628 & 0.742 & 0.991 & 0.994 & 0.998\\
0.01& 0.596 & 0.754 & 0.873 & 0.987 & 0.998 & 0.999\\ \hline
\end{tabular*}
\end{table*}

\section{Simulation study}\label{sec6}

We have indicated in the introduction that the original test for a
constant volatility from the noise-free model loses its asymptotic
properties in the
presence of noise. Unsurprisingly, for a smaller
variance of the noise variables, the data look more like observations
from a continuous semimartingale and thus the test statistics behaves
roughly in the same way as before, provided that the sample size is not
too large. On the other hand, for a large variance of the error terms
these are dominating, and thus the whole
procedure breaks down even for small sample sizes. The same problem
arises if the variance of the error is small but the sample size is
large (see the discussion in the \hyperref[sec1]{Introduction}).
We start with a further example simulating the level of the bootstrap
test proposed by Dette and
Podolskij~\cite{dp} for a parametric hypothesis, assessing its quality
for various sample sizes $n$ and different variances $\omega^2$.

Precisely, we have used that test for testing the hypothesis $H_0\dvt
\sigma^2(t,x) = \theta x^2$,
where $b(t,x)=0.1x$.
The results are obtained from 1000 simulation runs and 500 bootstrap
replications and displayed in Table~\ref{tab2} for various sample sizes and
standard deviations $\omega$ of the noise process. We observe that for
$n=256$ and a (small) standard deviation of $\omega=0.001$ the test
does roughly keep its asymptotic level, whereas it cannot be used at
all when the variance becomes larger. Moreover, even if the variance is
small but the sample size is increased, the test does not keep its
pre-assigned level
(see the results for $\omega=0.001$ and $n=1024$ in Table~\ref{tab2}). Thus, in
practice the application of testing procedures addressing the problem
of microstructure noise is strictly recommended.

In the following section, we illustrate the finite sample properties of
a bootstrap version of the Kolmogorov--Smirnov test based on the
processes investigated in Sections~\ref{sec4} and~\ref{sec5}. Since
the stochastic order of $|\Delta_i^n Z|$ is basically determined by
the maximum of $n^{-\fraca12}$ and $\omega$
(which are the orders of $|\Delta_i^n X|$ and $|\Delta_i^n U|$,
respectively), we kept $n \omega^2 = 0.1024$ fixed in order to have comparable
results for different sample sizes $n$. The regularisation parameters
$\kappa$ and $\rho$ were set to be $1/2$ each.
All simulation results presented in the following paragraphs are based
on 1000 simulation runs and 500 bootstrap replications (if the bootstrap
is applied to estimate critical values).

For all testing problems discussed below, we have not used exactly the
statistics $\hat N_t$ and $\hat M_t$, but related versions
accounting for finite sample adjustments. Following Jacod \textit{et
al.}~\cite{jlmpv}, where it has been shown that finite sample corrections
improve the behaviour of the estimate $\hat B_t^0$ (and presumably of
$\hat C$ as well) substantially, we have replaced the quantities $\psi
_i$ and $\Phi_{ij}$
in (\ref{psi}) by
certain numbers $\psi_i^n$ and $\Phi^n_{ij}$, which constitute the
``true'' quantities for finite samples, but are replaced by their limits
$\psi_i$ and $\Phi_{ij}$ in the asymptotics. See Jacod \textit{et al.}
\cite{jlmpv} for details.

\subsection{Testing for homoscedasticity}\label{sec6.1}

In the problem of testing for homoscedasticity the limiting process
$A(t)_{t \in[0,1]}$ has an extremely simple form, when the null
hypothesis of
a constant volatility holds. In fact, the finite dimensional
distributions of the process $(A(t))_{t \in[0,1]}$ coincide with those
of a rescaled Brownian bridge, thus $(A_n(t)/\hat s_t)_{t \in[0,1]}$
converges weakly to $(B_t)_{t \in[0,1]}$. We have investigated the
properties of the Kolmogorov--Smirnov test for different sample sizes
$n$, where the noise satisfies $U \sim\mathcal N(0, \omega^2)$ and the
drift function is again given by $b(t,x) = 0.1 x$. A similar test can
be constructed using Theorem~\ref{absol}, but the corresponding results
are omitted for the sake of brevity as the rate of convergence in this
case becomes worse.

\begin{table*}[b]
\tabcolsep=0pt
\tablewidth=270pt
\caption{Simulated nominal level of the test, which rejects the null
hypothesis of homoscedasticity for a large value of
$\sup|A_n(t)/\hat s_t|$, using the critical values from the asymptotic
theory. The variance
of the noise process is defined by $n \omega^2= 0.1024$}
\label{tab3}\begin{tabular*}{270pt}{@{\extracolsep{\fill}}d{5.0}d{1.3}d{1.3}d{1.3}@{}}
\hline
\multicolumn{1}{@{}l}{{$n$}/{$\alpha$}} & 0.025 & 0.05 & 0.1 \\ \hline
256 & 0.008 & 0.022 & 0.058 \\
1024 & 0.007 & 0.023 & 0.062 \\
4096 & 0.013 & 0.029 & 0.079 \\
16\,384 & 0.017 & 0.038 & 0.077 \\ \hline
\end{tabular*}
\end{table*}

In Table~\ref{tab3}, we present the simulated level of the Kolmogorov--Smirnov
test using the critical values from the asymptotic distribution.
It can be seen that the asymptotic level of the test is slightly
underestimated. This effect
becomes less visible for a larger sample size, but even then it is
still apparent. Note that these findings are in line with previous
simulations on noisy observations and it is likely that they are due to
the fact the rate of convergence for most testing problems is only
$n^{-\fraca14}$.

\subsection{Testing general hypotheses}\label{sec6.2}

For a general null hypothesis in (\ref{h0}), the distribution of the
limiting process $(A(t))_{t \in[0,1]}$ depends on the path of the
underlying semimartingale $(X_t)_{t \in[0,1]}$ and on the volatility
$(\sigma_t)_{t \in[0,1]}$, and thus we cannot use it directly for the
calculation
of critical values. For this reason, we propose the application of the
parametric bootstrap in order to obtain simulated
critical values. First, we compute the global estimators $\hat
\omega^2$ and $\hat\theta= \hat D^{-1} \hat C$ as well as each
$n^{\fraca14} \hat N_t$ and $\hat s^2_t$ from the observed data. Under
the null hypothesis $N_t$ equals zero,
and thus it is intuitively clear that the null hypothesis has to be
rejected for large values of the standardised Kolmogorov--Smirnov statistic
$Y_n = \sup_{t \in[0,1]} |n^{\fraca14} \hat N_t/\hat s_t|$.

\begin{table*}[b]
\tabcolsep=0pt
\tablewidth=300pt
\caption{Simulated level of the bootstrap test based on the
standardised Kolmogorov--Smirnov functional of $(\hat N_t)$ for various
hypotheses. The variance
of the noise process is defined by $n \omega^2= 0.1024$}
\label{tab4}\begin{tabular*}{300pt}{@{\extracolsep{\fill
}}d{4.0}d{1.3}d{1.3}d{1.3}d{1.3}d{1.3}d{1.3}@{}}
\hline
\multicolumn{1}{@{}l}{$\sigma_1^2(t,x)$} & \multicolumn{3}{l}{1} &
\multicolumn{3}{l@{}}{$x^2$} \\[-5pt]
& \multicolumn{3}{l}{\hrulefill} & \multicolumn{3}{l@{}}{\hrulefill}\\
\multicolumn{1}{@{}l}{{$n$}/{$\alpha$}} & 0.025 & 0.05 & 0.1 & 0.025 &
0.05 & 0.1 \\
\hline
256 & 0.019 & 0.046 & 0.113 & 0.03 & 0.066 & 0.118 \\
1024 & 0.02 & 0.049 & 0.099 & 0.034 & 0.07 & 0.119 \\
4096 & 0.021 & 0.04 & 0.072 & 0.022 & 0.048 & 0.090 \\ \hline
\end{tabular*}
\end{table*}

In a second step we generate bootstrap data $Z^{*(j)}_{\fraca1n} =
X^{*(j)}_{\fraca1n} + U^{*(j)}_{\fraca1n}$, where the $X^{*(j)}_{\fraca
in}$ are realisations of the process in (\ref{2.3}) with $b_s \equiv0$
and $\sigma^2_s = \sigma^2(s,X_s)
= \sum_{k=1}^{d} \hat\theta_k \sigma_k^2(s,X_s)$ (corresponding to the
null hypothesis) and each $U^{*(j)}_{\fraca in}$ is normally distributed
with mean zero and variance $\hat\omega^2$. Using these data, we
calculate the corresponding bootstrap statistics $Y_n^{*(j)}$ and use
these to
compute the quantiles of the bootstrap distribution. Finally, the null
hypothesis is rejected if $Y_n$ is larger than its $(1- \alpha)$-quantile.

In order to investigate the approximation of the nominal level we
consider the hypothesis of constant volatility and the hypothesis
$H_0\dvt \sigma^2(t,x)=\theta x^2$. The data is generated under the null
hypothesis with drift function
$b(t,x)=0.1x$ and the rejection probabilities are depicted in Table~\ref{tab4}.
These results show that the bootstrap approximation works
well even for a small $n$. In particular, we see that in the case of
homoscedasticity the
exact asymptotic test using the weak convergence of $Y_n$ to the
supremum of a standard Brownian bridge is outperformed (compare with
Table~\ref{tab3}).
In the case of testing, the parametric hypothesis $H_0\dvt \sigma
^2(t,x)=x^2$ we observe a slight overestimation of the nominal level by
the bootstrap test.

As an example for testing the hypothesis $\bar H_0$, we have chosen
$\sigma(t,x) = \theta|x|$ and investigated the properties of the
analogues of $Y_n$ and
$Y_n^{*(j)}$ from above, where we have replaced $n^{\fraca14} \hat N_t$
and $\hat s_t$ by $n^{\fraca14 - \fraca\delta2} \hat M_t$ and $\hat
r_t$, respectively. In this case, we chose $\delta= \frac14$,
corresponding to $l_n = \mathrm{O}(n^{-\fraca34})$ and a rate of convergence
$n^{-\fraca
18}$. Note that in this particular situation there is no need for
stating the hypothesis in terms of $\bar H_0$ as it is equivalent to
$\sigma^2(t,x) = \theta|x|^2$, but nevertheless it gives a reasonable
impression on how well the bootstrap approximation works for testing
hypotheses of the form $\bar H_0$.

\begin{table*}
\tabcolsep=0pt
\tablewidth=200pt
\caption{Simulated level of the bootstrap test based on the
standardised Kolmogorov--Smirnov functional of $(\hat M_t)$ for
$\sigma(t,x) = \theta|x|$. The variance
of the noise process is defined by $n \omega^2= 0.1024$}
\label{tab5}\begin{tabular*}{200pt}{@{\extracolsep{\fill
}}d{4.0}d{1.3}d{1.3}d{1.3}@{}} \hline
\multicolumn{1}{@{}l}{{$n$}/{$\alpha$}} & 0.025 & 0.05 & 0.1 \\
\hline
256 & 0.040 & 0.076 & 0.136 \\
1024 & 0.032 & 0.057 & 0.119 \\ \hline
\end{tabular*}
\end{table*}

\begin{table*}[b]
\tabcolsep=0pt
\caption{Simulated rejection probabilities of the bootstrap test based
on the standardised Kolmogorov--Smirnov functional of $(\hat
N_t)$ for various alternatives. The data is simulated with $\sigma
^2(t,x) = \theta|x|^2$ and the variance of the noise process
is defined by $n \omega^2=0.1024$}
\label{tab6}\begin{tabular*}{\textwidth}{@{\extracolsep{\fill
}}d{4.0}d{1.3}d{1.3}d{1.3}d{1.3}d{1.3}d{1.3}d{1.3}d{1.3}d{1.3}@{}}
\hline
\multicolumn{1}{@{}l}{alt} & \multicolumn{3}{l}{1} & \multicolumn
{3}{l}{$1+|x|$} & \multicolumn{3}{l@{}}{Heston} \\[-5pt]
& \multicolumn{3}{l}{\hrulefill} & \multicolumn{3}{l}{\hrulefill} &
\multicolumn{3}{l@{}}{\hrulefill}
\\
\multicolumn{1}{@{}l}{{$n$}/{$\alpha$}} & 0.025 & 0.05 & 0.1 & 0.025 &
0.05 & 0.1 & 0.025 & 0.05 & 0.1 \\ \hline
256 & 0.057 & 0.128 & 0.237 & 0.073 & 0.152 & 0.263 & 0.722 & 0.870 &
0.941 \\
1024 & 0.170 & 0.230 & 0.329 & 0.224 & 0.326 & 0.465 & 0.975 & 0.980 &
0.985 \\ \hline
\end{tabular*}
\end{table*}

We observe from the results in Table~\ref{tab5} that even though the rate of
convergence in Theorem~\ref{absol} is worse than in Theorem~\ref{tigh},
there is no
substantial difference in the approximation of the nominal level by the
bootstrap test for both types of hypotheses: The nominal level is
slightly overestimated, but in general
the parametric bootstrap yields a satisfactory and reliable
approximation of the nominal level.

Finally, Table~\ref{tab6} contains the rejection probabilities of the bootstrap
test under the alternative. The null hypothesis is given by $H_0\dvt \sigma
^2(t,x) = \theta|x|^2$, and we discuss two local volatility
alternatives, namely $\sigma^2(t,x) = 1$ and $\sigma^2(t,x) = 1 + |x|$,
and one alternative
coming from a stochastic volatility model is considered. For this case,
we chose the Heston model, that is,
\begin{eqnarray*}
&&X_t = X_0 + \int_0^t (\mu-\nu_s/2) \,\mathrm{d}s + \int_0^t \sigma_s \,\mathrm{d}W_t
\\
&& \quad \mbox{with }  \nu_t = \nu_0 + \delta\int_0^t (\alpha-\nu_s) \,\mathrm{d}s +
\gamma\int_0^1 \nu_s^{1/2}\,\mathrm{d}B_s,
\end{eqnarray*}
where
$\nu_t=\sigma_t^2$ and $\Corr(W, B)=\eta$ and the parameters were
chosen as $\mu=0.05/252, \delta=5/252,\alpha=0.04/252,
\gamma=0.05/252$ and $\rho=-0.5$.

We observe from the results depicted in Table~\ref{tab6} that the bootstrap test
indicates in all cases that the null hypothesis is not satisfied. It is
also remarkable that it is more difficult to detect the local
volatility alternatives than the one coming from the Heston model. In
the latter case, the rejection probabilities are extremely large even
for a small sample size, contrary to the first two situations.

\begin{appendix}\label{appm}
\section*{\texorpdfstring{Appendix: Proof of Theorem \protect\ref{tigh}}{Appendix: Proof of Theorem 1}}\label{sec7}

We will only prove the Theorem~\ref{tigh}, as similar methods show
Theorem~\ref{absol} as well. We start with a typical localisation
argument, which allows us to assume that several quantities are
bounded. Recall first that $a$ and $\sigma$ are locally bounded by
assumption, from which is follows that
$X$ is locally bounded as well. Thus we can conclude along the lines of
Jacod~\cite{jac08} that we may assume without loss of generality that
each of these processes is actually bounded. Since further each $\sigma
^2_i$ is continuous and because $U$ has a compact support, we may
conclude that both $(s,X_t)$ and $(s,\hat X_{\fraca kn})$ (for arbitrary
$s$, $t$, $k$ and $n$) are living on a compact set, and thus
$\sigma^2_i(s,X_t)$ and $\sigma^2_i(s, \hat X_{\fraca kn})$ are also
bounded, the latter one uniformly in $n$. Similar results hold for the
first two derivatives of $\sigma_i^2$ as well as for any of the
functions $\bar\sigma_i$. Constants are denoted by $K$ throughout this
section.

The proof of Theorem~\ref{tigh} 
is based on several preliminary results, and we start with two results
determining the rate of convergence of the quantities $\hat{B}_t^i -
B_t^i$ and $\hat{D}_{ij} - D_{ij}$ defined in (\ref{2.21a}) and (\ref
{dt}), respectively. The following result ensures that the
(conditional) variance in a limit theorem for $\hat{N}_t - N_t$ will
not depend on $\hat{B}_t^i$ and $\hat{D}_{ij}$, since the rate of
convergence is $n^{-\fraca14}$. Thus, we will focus in the following on
the behavior of $\hat{C}_i$ and $\hat{B}_t^0$.

\begin{theorem} \label{klord}
Under the assumptions from Section~\ref{assum} we have
%
\begin{eqnarray} \label{est1}
\hat{B}_t^i - B_t^i &=& \mathrm{o}_p(n^{-\fraca14}) \qquad\mbox{for }  i=1, \ldots
, d, \nonumber
\\[-8pt]
\\[-8pt]\hat{D}_{ij} - D_{ij} &=& \mathrm{o}_p(n^{-\fraca14}) \qquad\mbox{for }
 i,j=1, \ldots, d,
\nonumber
\end{eqnarray}
where the first result holds uniformly with respect to $t \in[0,1]$.
\end{theorem}

\begin{pf} For a proof of the first estimate,
we use for a fixed index $i$ the decomposition
\begin{eqnarray*}
\hat{B}_t^i - B_t^i &=& {\frac1n} \sum_{k=1}^{\lfloor nt \rfloor
-m_n}  \biggl(\sigma_i^2\biggl ({\frac kn}, \hat{X}_{\fraca kn}\biggr) - \sigma_i^2
\biggl({\frac
kn}, X_{\fraca kn}\biggr)  \biggr)\\
&&{} +  \Biggl({\frac1n} \sum_{k=1}^{\lfloor nt
\rfloor-m_n} \sigma_i^2 \biggl({\frac kn}, X_{\fraca kn}\biggr) - \int_0^t \sigma_i^2
(s, X_s) \,\mathrm{d}s \Biggr).
\end{eqnarray*}
Regarding the first term in this sum, note that
\begin{eqnarray*}
\hat{X}_{\fraca kn} - X_{\fraca kn} = \frac{1}{m_n} \sum_{j=1}^{m_n}
\biggl(U_{\fracb{k+j}{n}} + \int_{\fraca kn}^{\fracb{k+j}{n}} \sigma_u \,\mathrm{d}W_u
 \biggr) + \mathrm{O}_p(n^{-\fraca12})
\end{eqnarray*}
and thus $ \hat{X}_{\fraca kn} - X_{\fraca kn} = \mathrm{O}_p(n^{-\fraca14}).$ A
Taylor expansion and boundedness of the second derivative of the
function $\sigma^2$ give
\begin{eqnarray*}
\frac1n \sum_{k=1}^{\lfloor nt \rfloor-m_n}  \biggl(\sigma_i^2 \biggl({\frac
kn}, \hat{X}_{\fraca kn}\biggr) - \sigma_i^2 \biggl({\frac kn}, X_{\fraca kn}\biggr)  \biggr)
= \frac1n \sum_{k=1}^{\lfloor nt \rfloor-m_n} A_{k,n} + \mathrm{O}_p(n^{-\fraca12})
\end{eqnarray*}
with
\begin{eqnarray*}
A_{k,n} = \frac{1}{m_n} \sum_{j=1}^{m_n} \frac{\partial}{\partial y}
\sigma_i^2 \biggl({\frac kn}, X_{\fraca kn}\biggr)   \biggl(U_{\fracb{k+j}{n}} + \int
_{\fraca kn}^{\fracb{k+j}{n}} \sigma_s \,\mathrm{d}W_s \biggr).
\end{eqnarray*}
However, we have $E[A_{k,n} A_{l,n}] = \mathrm{O}(n^{-\fraca12})$ for arbitrary
$k$ and $l$ as well as $E[A_{k,n} A_{k+l,n}] = 0$ for $l \geq m_n$ by
conditioning on $\mathcal F_{\fracb{k+l}{n}}$. This yields
\begin{eqnarray*}
E \Biggl[ \Biggl({\frac1n} \sum_{k=1}^{\lfloor nt \rfloor-m_n} A_{k,n}
\Biggr)^2 \Biggr] = \frac{1}{n^2} \sum_{k=m_n}^{\lfloor nt \rfloor-2m_n} \sum
_{l=-m_n}^{m_n} E[A_{k,n} A_{k+l,n}] + \mathrm{O}\biggl(\frac{m_n}{n^2}\biggr) = \mathrm{O}\biggl(\frac1n\biggr),
\end{eqnarray*}
which is small enough. For the second term in the decomposition of $\hat
{B}_t^i - B_t^i$ it holds that
\begin{eqnarray*}
&& {\frac1n} \sum_{k=1}^{\lfloor nt \rfloor-m_n} \sigma_i^2 \biggl({\frac
kn}, X_{\fraca kn}\biggr) - \int_0^t \sigma_i^2 (s, X_s) \,\mathrm{d}s \\
&& \quad = \sum
_{k=1}^{\lfloor nt \rfloor} \int_{\fracb{k-1}{n}}^{\fraca kn}  \biggl(\sigma
_i^2 \biggl({\frac{k-1}{n}}, X_{\fracb{k-1}{n}}\biggr) - \sigma_i^2 \bigl(s,
X_{\fracb
{k-1}{n}}\bigr)\\
&&\hphantom{\quad = \sum
_{k=1}^{\lfloor nt \rfloor} \int_{\fracb{k-1}{n}}^{\fraca kn}  \biggl(}{} + \sigma_i^2 \bigl(s, X_{\fracb{k-1}{n}}\bigr) - \sigma_i^2 (s, X_s)
 \biggr) \,\mathrm{d}s + \mathrm{O}_p(n^{-\fraca12}).
\end{eqnarray*}
By differentiability in both components and from a similar expansion as
above the claim follows. The result on $\hat{D}_{ij} - D_{ij}$ can be
shown in the same way.
\end{pf}

The following result specifies the convergence of the finite
dimensional distributions of the processes, which are used for the
construction of $\{\hat N_t \}_{t\in[0,1]}$. Below we use the notation
$G_n \stab G$ to indicate stable convergence of a sequence of random
variables $(G_n)$ to a limiting variable $G$, which is defined on an
appropriate extension of the original probability space. For details on
stable convergence see Jacod and Shiryaev
\cite{jacshir}.

\begin{theorem} \label{fin}
Define for any fixed $t_1, \ldots, t_k \in[0,1]$ the matrix
$\Sigma_{t_1, \ldots, t_k} (s,X_s) = \break\gamma_s^2  \ell(s,X_s) \ell^T
(s, X_s) $ where $\ell(s,X_s) = (
1_{[0,t_1]}(s), \ldots, 1_{[0,t_k]}(s) ,h^T(s,X_s))^T $. Then we have
\begin{eqnarray*} n^{\fraca14}  ( \hat{B}_{t_1}^{0} - B_{t_1}^0 ,
\ldots,
\hat{B}_{t_k}^{0} - B_{t_k}^0, \hat{C}_1 - C_1 , \ldots,\hat{C}_d -
C_d  )^T
\stab\int_0^1 \Sigma^{\fraca12}_{t_1, \ldots, t_k} (s,X_s) \,\mathrm{d}W'_s,
\end{eqnarray*}
where $W'$ is another Brownian motion, which is independent of the
$\sigma$-algebra $\mathcal F$.
\end{theorem}

\begin{pf}
Since $\omega^2 - \hat{\omega}^2_n = \mathrm{O}_p(n^{-\fraca12})$, one obtains
\begin{eqnarray*}
 \hat{C}_i &=& \frac{1}{n} \sum_{k=1}^{n-m_n} \sigma_i
^2\biggl({\frac kn},X_{\fraca kn}\biggr)  \hat\sigma^2_{\fraca kn} +
\frac{1}{n} \sum_{k=1}^{n-m_n}  \biggl(\sigma_i ^2\biggl({\frac kn},\hat
{X}_{\fraca kn}\biggr) - \sigma_i ^2\biggl({\frac kn},X_{\fraca
kn}\biggr) \biggr)  \hat\sigma^2_{\fraca kn} \\
&&{}+ \mathrm{O}_p(n^{-\fraca12}).
\end{eqnarray*}
From similar arguments as given in the proof of Theorem~\ref{klord} we
find that the second term is of order $\mathrm{o}_p(n^{- \fraca14})$ and thus
asymptotically negligible as well. Therefore, we are left to focus on
$F_{in} = \frac{1}{n} \sum_{k=1}^{n-m_n} \sigma_i ^2({\frac
kn},X_{\fraca kn})  \hat\sigma^2_{\fraca kn}$. Due to the dependence
structure of the summands in $F_{in}$ it will be convenient to use a
``small-blocks--big-blocks''-technique as in Jacod \textit{et al.} \cite
{jlmpv} in order to prove Theorem~\ref{fin}. To this end, we choose an
integer $p$, which eventually goes to infinity, and partition the $n$
observations into several subsets: We define $b_k(p) = k(p+1)m_n$ and
$c_k(p) = k(p+1)m_n+ pm_n$ and denote by $j_n(p)$ the largest integer
$k$ such that $c_k(p) \leq n-m_n$ holds. Moreover, we use the notation
$i_n(p) = (j_n(p)+1) pm_n$, and introduce for each $0 \leq k \leq
j_n(p)$ and any $p$ the following random
variables:
\begin{eqnarray*}
G(k,p)_1^n &=& \frac{1}{n} \sigma_i^2\biggl(\frac{b_k(p)}{n}, X_{\fraca
{b_k(p)}{n}}\biggr) \sum_{j=b_k(p)}^{c_k(p)-1} \hat\sigma^2_{\fraca kn},
\\
G(k,p)_2^n &=& \frac{1}{n} \sigma_i^2\biggl(\frac{c_k(p)}{n}, X_{\fraca
{c_k(p)}{n}}\biggr) \sum_{j=c_k(p)}^{b_{k+1}(p)-1} \hat\sigma^2_{\fraca kn}.
\end{eqnarray*}
The remainder terms from $i_n(p)$ to $n-m_n$ are gathered in some
$G(p)_3^n$. Note that each of these quantities depends on $i$, although
it does not appear in the notation.

The main intuition behind these quantities is that the terms
$G(k,p)_1^n$ are defined on non-overlapping intervals, which means that
the intervals on which each $\overline{Z}^n_j$ within $G(k,p)_1^n$
lives are disjoint from those of any $\overline{Z}^n_j$ within any
other $G(l,p)_1^n$. This is sufficient to ensure some type of
conditional independence, which will be used in order to prove Theorem
\ref{fin}. The variables $G(k,p)_2^n$ and $G(p)_3^n$ are filling the
gaps between $G(k,p)_1^n$ and $G(l,p)_1^n$ and can be shown to be
asymptotically negligible.

An important tool will be the following decomposition of $|\overline
{Z}^n_j|^2$. We set
\begin{eqnarray*} V^j_s = \int_{\fraca jn}^{\fraca jn+s} g_n\biggl(u - \frac
jn\biggr)  a_u \,\mathrm{d}u + \int_{\fraca jn}^{\fraca jn+s} g_n\biggl(u - \frac jn\biggr)  \sigma_u
\,\mathrm{d}W_u,
\end{eqnarray*}
and obtain by an application of Ito's formula
\begin{eqnarray*} \label{ito}
|\overline{Z}^n_j|^2 &=& |\overline
{X}^n_j|^2 + |\overline{U}^n_j|^2 + 2\overline{X}^n_j  \overline{U}^n_j
\\
&=& 2 \int_{\fraca jn}^{\fracb{j + m_n}{n}} V^j_s  g_n\biggl(s -
\frac jn\biggr)  a_s \,\mathrm{d}s + 2 \int_{\fraca jn}^{\fracb{j + m_n}{n}} V^j_s  g_n\biggl(s -
\frac jn\biggr)  \sigma_s \,\mathrm{d}W_s\\
&&{} + \int_{\fraca jn}^{\fracb{j +
m_n}{n}} g_n^2\biggl(s - \frac jn\biggr)  \sigma^2_s \,\mathrm{d}s + |\overline
{U}^n_j|^2 + 2
\overline{U}^n_j \int_{\fraca jn}^{\fracb{j + m_n}{n}}  g_n\biggl(s - \frac
jn\biggr)  a_s \,\mathrm{d}s\\
&&{} + 2 \overline{U}^n_j \int_{\fraca jn}^{\fracb{j +
m_n}{n}}  g_n\biggl(s - \frac jn\biggr)  \sigma_s \,\mathrm{d}W_s\\
 &=& \sum_{l=1}^{6} D(j)_l^n,
\end{eqnarray*}
where the last identity defines the quantities $D(j)^n_l$ in an
obvious manner.

For $b_k(p) \leq j < c_k(p)$ we introduce approximations for the
quantities $D(j)_2^n$ and $D(j)_6^n$, namely
\begin{eqnarray*}
\tilde D(k,j,p)_2^n &=& 2 \sigma^2_{\fraca{b_k(p)}{n}} \int_{\fraca
jn}^{\fracb{j + m_n}{n}}  \biggl( \int_{\fraca jn}^{\fraca jn+s} g_n\biggl(u -
\frac jn\biggr) \,\mathrm{d}W_u  \biggr)
g_n\biggl(s - \frac jn\biggr) \,\mathrm{d}W_s, \\
\tilde D(k,j,p)_6^n &=& 2 \sigma_{\fraca{b_k(p)}{n}}  \overline{U}^n_j
\int_{\fraca jn}^{\fracb{j + m_n}{n}}  g_n\biggl(s - \frac jn\biggr)\,\mathrm{d}W_s.
\end{eqnarray*}
Additionally, we set $H(k,p)^n = \sigma_i^2(\frac{b_k(p)}{n}, X_{\fraca
{b_k(p)}{n}})  Y(k,p)^n$, where
%
\begin{equation}
 \hspace*{-15pt}\label{defY} Y(k,p)^n = \frac{1}{\kappa\psi_2} n^{-\fraca12}
\sum_{j=b_k(p)}^{c_k(p)-1}
 \biggl\{ \tilde D(k,j,p)_2^n + \tilde D(k,j,p)_6^n +  \biggl(D(j)_4^n -
n^{-\fraca12} \frac{\psi_1}{\kappa} \omega^2 \biggr)  \biggr\}.
\end{equation}
Finally, we
define
\begin{eqnarray*} \chi(p)_k^n = E\Bigl[\Bigl(\sup_{s, t \in[b_k(p)/n, c_k(p)/n]}
|a_s - a_t| + |\sigma_s - \sigma_t|\Bigr)^2
 \big| \mathcal F_{\fraca{b_k(p)}{n}}\Bigr]^{\fraca12} .
\end{eqnarray*}

The main part of the proof of Theorem~\ref{tigh} are two auxiliary
results which specify the asymptotic properties of $F_{in}$.

\begin{lemma} \label{lem1}
We have
\begin{eqnarray*} \lim_{p \rightarrow\infty} \limsup_{n \rightarrow
\infty} n^{\fraca14} \Biggl \{ \Biggl (\sum_{k=0}^{j_n(p)} \bigl(G(k,p)_1^n + G(k,p)_2^n\bigr)
+ G(p)_3^n - C_i  \Biggr) - \sum_{k=0}^{j_n(p)} H(k,p)  \Biggr\} = 0.
\end{eqnarray*}
\end{lemma}
\begin{pf}   The proof goes through a rather
large number of steps and makes extensive use of the decomposition in
(\ref{ito}). We will show first that the influence of the random
variables $D(j)_1^n$ and $D(j)_5^n$ within $G(k,p)_1^n$ (and
analogously for $G(k,p)_2^n$ and $G(p)_3^n$) is asymptotically
negligible, that is
%
\begin{eqnarray}
\label{1.1}  \lim_{p \rightarrow\infty} \limsup_{n \rightarrow
\infty} n^{-\fraca14} \sum_{k=0}^{j_n(p)}
\sigma_i^2\biggl(\frac{b_k(p)}{n}, X_{\fraca{b_k(p)}{n}}\biggr) \sum
_{j=b_k(p)}^{c_k(p)-1} \bigl(D(j)_1^n + D(j)_5^n\bigr) = 0.
\end{eqnarray}
For a proof of (\ref{1.1}), assume without loss of generality that
$b_k(p) \leq j < c_k(p)$. One obtains
\begin{eqnarray*}
D(j)_1^n &=& 2 a_{\fraca{b_k(p)}{n}} \int_{\fraca jn}^{\fracb{j + m_n}{n}}
\biggl (
\int_{\fraca jn}^{\fraca jn +s}  g_n\biggl(u - \frac jn\biggr)  \sigma_u \,\mathrm{d}W_u  \biggr)
 g_n\biggl(s - \frac jn\biggr) \,\mathrm{d}s \\
 &&{}+ 2 \int_{\fraca jn}^{\fracb{j + m_n}{n}}
 \biggl( \int_{\fraca jn}^{\fraca jn +s}  g_n\biggl(u - \frac jn\biggr)  \sigma_u \,\mathrm{d}W_u
 \biggr)  g_n\biggl(s - \frac jn\biggr) \bigl(a_s - a_{\fraca{b_k(p)}{n}}\bigr) \,\mathrm{d}s \\
 &&{}+ \mathrm{O}_p\biggl(\frac1n\biggr)
\end{eqnarray*}
and from the martingale property of a stochastic integral with respect
to Brownian motion and the Cauchy--Schwarz inequality we derive that
$|E[D(j)_1^n  |\mathcal F_{\fraca{b_k(p)}{n}}]| \leq K  n^{-\fraca34}  \chi(p)_k^n.$ Thus, with the notation $\delta(k,p)_1^n = \sum
_{j=b_k(p)}^{c_k(p)-1} D(j)_1^n$ we conclude
\begin{eqnarray*}
\bigl|E\bigl[\delta(k,p)_1^n  |\mathcal
F_{\fraca{b_k(p)}{n}}\bigr]\bigr| \leq K  p  n^{-\fraca14}  \chi(p)_k^n \quad
\mbox{and} \quad E\bigl[(\delta(k,p)_1^n)^2  | \mathcal F_{\fraca
{b_k(p)}{n}}\bigr] \leq K  p^2  n^{-\fraca12},
\end{eqnarray*}
and for $k > l$ it follows
\begin{eqnarray*}  &&\biggl| E \biggl\{\sigma_i^2\biggl(\frac{b_k(p)}{n},
X_{\fraca
{b_k(p)}{n}}\biggr)  \sigma_i^2\biggl(\frac{b_l(p)}{n}, X_{\fraca{b_l(p)}{n}}\biggr)  \delta(l,p)_1^n  E \bigl[\delta(k,p)_1^n
| \mathcal F_{\fraca{b_k(p)}{n}} \bigr]  \biggr\}
 \biggr|\\
 && \quad  \leq K  p^2  n^{-\fraca12}  E[\chi(p)_k^n].
\end{eqnarray*}
Since $j_n(p)$ is of order ${n^{1/2}/p}$, we obtain
\begin{eqnarray*} &&E \Biggl[\Biggl (n^{-\fraca14} \sum_{k=0}^{j_n(p)} \sigma
_i^2\biggl(\frac{b_k(p)}{n}, X_{\fraca{b_k(p)}{n}}\biggr) \sum_{j=b_k(p)}^{c_k(p)-1}
D(j)_1^n \Biggr)^2 \Biggr] \\
&& \quad \leq K  \Biggl( p  n^{-\fraca12} + \sum
_{k>l}^{j_n(p)} p^2  n^{-1}  E[\chi(p)_k^n]  \Biggr).
\end{eqnarray*}
From Lemma 5.4. in Jacod \textit{et al.}~\cite{jlmpv} it follows that $
\lim_{n \rightarrow\infty} n^{-\fraca12}
\sum_{k=1}^{j_n(p)} E[\chi(p)_k^n] = 0 $ for any $p$, which gives that
the first term in the sum (\ref{1.1}) converges to 0. The second term
in (\ref{1.1}) converges to zero from the independence of $X$ and $U$
and a standard martingale argument.

The next step is devoted to the analysis of the term $D(j)_2^n$. We prove
%
\begin{equation}
 \label{1.2.1} \lim_{p \rightarrow\infty} \limsup_{n \rightarrow
\infty} n^{-\fraca14} \sum_{k=0}^{j_n(p)} \sigma_i^2\biggl(\frac{b_k(p)}{n},
X_{\fraca{b_k(p)}{n}}\biggr) \sum_{j=b_k(p)}^{c_k(p)-1}  \bigl(D(j)_2^n - \tilde
D(k,j,p)_2^n \bigr) = 0
\end{equation}
as well as
%
\begin{eqnarray} \label{1.2.2}  \lim_{p \rightarrow\infty} \limsup_{n
\rightarrow\infty} n^{-\fraca14} \sum_{k=0}^{j_n(p)} \sigma_i^2\biggl(\frac
{c_k(p)}{n}, X_{\fraca{c_k(p)}{n}}\biggr) \sum_{j=c_k(p)}^{b_{k+1}(p)-1}
D(j)_2^n &=& 0, \\
\label{1.2.2b}  \lim_{p \rightarrow\infty} \limsup
_{n \rightarrow\infty} n^{-\fraca14}  \sigma_i^2\biggl(\frac{i_n(p)}{n},
X_{\fraca{i_n(p)}{n}}\biggr) \sum_{j=i_n(p)}^{n-m_n} D(j)_2^n &=& 0.
\end{eqnarray}
Set $b_k(p) \leq j < c_k(p)$ again. A martingale argument as before
allows us to focus on
\begin{eqnarray*}
D''(j)_2^n = 2 \int_{\fraca jn}^{\fracb{j + m_n}{n}}  \biggl( \int_{\fraca
jn}^{\fraca jn +s}  g_n\biggl(u - \frac jn\biggr)  \sigma_u \,\mathrm{d}W_u  \biggr)  g_n\biggl(s -
\frac jn\biggr)  \sigma_s \,\mathrm{d}W_s
\end{eqnarray*}
only. We have $E[D''(j)_2^n  | \mathcal F_{\fraca{b_k(p)}{n}}] = 0$ and
$E[|D''(j)_2^n D''(l)_2^n|  | \mathcal F_{\fraca{b_k(p)}{n}}] \leq K  n^{-1}$,
thus (\ref{1.2.2b}) follows easily. For (\ref{1.2.2}), note that
$E[(\sum_{j=c_k(p)}^{b_{k+1}(p)-1} D''(j)_2^n)^2] \leq K, $ which gives
(recall the definition of $j_n(p), b_k(p)$ and $c_k(p)$)
\begin{eqnarray*}
n^{-\fraca12} \sum_{k=0}^{j_n(p)} E\Biggl [\sigma_i^4\biggl(\frac{c_k(p)}{n},
X_{\fraca{c_k(p)}{n}}\biggr)  \Biggl(\sum_{j=c_k(p)}^{b_{k+1}(p)-1} D''(j)_2^n
\Biggr)^2 \Biggr]
\leq K  n^{-\fraca12}  \frac{n^{\fraca12}}{p} = K  \frac1p,
\end{eqnarray*}
converging to zero as $p$ tends to infinity. We are thus left to prove
\begin{eqnarray*}
\lim_{p \rightarrow\infty} \limsup_{n \rightarrow\infty} n^{-\fraca
14} \sum_{k=0}^{j_n(p)} \sigma_i^2\biggl(\frac{b_k(p)}{n}, X_{\fraca
{b_k(p)}{n}}\biggr) \sum_{j=b_k(p)}^{c_k(p)-1} \bigl (D''(j)_2^n - \tilde
D(k,j,p)_2^n \bigr) = 0.
\end{eqnarray*}
This time, we have $ E[D''(j)_2^n -\tilde D(k,j,p)_2^n  | \mathcal
F_{\fraca{b_k(p)}{n}}] = 0 $ and
\begin{eqnarray*} E\bigl[\bigl|\bigl(D''(j)_2^n- \tilde D(k,j,p)_2^n\bigr)  \bigl(D''(l)_2^n -
\tilde D(k,l,p)_2^n\bigr) \bigr|  | \mathcal F_{\fraca{b_k(p)}{n}}\bigr] \leq K  n^{-1}  (\chi(p)_k^n)^2.
\end{eqnarray*}
Thus,
\begin{eqnarray*}
&&E \Biggl[ \Biggl\{n^{-\fraca14} \sum_{k=0}^{j_n(p)} \sigma_i^2\biggl(\frac
{b_k(p)}{n}, X_{\fraca{b_k(p)}{n}}\biggr) \sum_{j=b_k(p)}^{c_k(p)-1}
\bigl(D(j)_2^n -
\tilde D(k,j,p)_2^n \bigr) \Biggr\}^2 \Biggr] \\
&& \quad \leq K p^2  n^{-\fraca12} \sum
_{k=0}^{j_n(p)} E [(\chi(p)_k^n)^2 ],
\end{eqnarray*}
and with a similar argument as in the proof of (\ref{1.1}) we are done.
Proving that $D(j)_6^n$ can be replaced by $\tilde D(k,j,p)_6^n$ works
analogously, thus we finish the proof of Lemma~\ref{lem1} showing
%
\begin{eqnarray} \label{1.3}
\hspace*{-25pt}&&\lim_{p \rightarrow\infty} \limsup_{n \rightarrow\infty}
n^{\fraca14}  \Biggl\{\frac{1}{\kappa\psi_2} n^{-\fraca12}  \Biggl(\sum
_{k=0}^{j_n(p)}  \Biggl( \sigma_i^2\biggl(\frac{b_k(p)}{n}, X_{\fraca{b_k(p)}{n}}\biggr)
\sum_{j=b_k(p)}^{c_k(p)-1} D(j)_3^n \nonumber
\\
\hspace*{-25pt}&&\hphantom{\lim_{p \rightarrow\infty} \limsup_{n \rightarrow\infty}
n^{\fraca14}  \Biggl\{\frac{1}{\kappa\psi_2} n^{-\fraca12}  \Biggl(\sum
_{k=0}^{j_n(p)}  \Biggl(}{}  + \sigma
_i^2\biggl(\frac{c_k(p)}{n}, X_{\fraca{c_k(p)}{n}}\biggr) \sum
_{j=c_k(p)}^{b_{k+1}(p)-1} D(j)_3^n
 \Biggr) \\
\hspace*{-25pt} &&\hphantom{\lim_{p \rightarrow\infty} \limsup_{n \rightarrow\infty}
n^{\fraca14}  \Biggl\{\frac{1}{\kappa\psi_2} n^{-\fraca12}  \Biggl(}{}+
\sum_{j=i_n(p)}^{n-m_n} D(j)_3^n  \Biggr)- C_i  \Biggr\} = 0.
\nonumber
\end{eqnarray}
We start with the following proposition:
%
\begin{eqnarray} \label{prop1}
&& \lim_{p \rightarrow\infty} \limsup_{n \rightarrow\infty} n^{\fraca
14}  \Biggl\{  \Biggl(\sum_{k=0}^{j_n(p)}  \biggl(
\int_{\fraca{b_k(p)}{n}}^{\fraca
{c_k(p)}{n}} \sigma_i^2\biggl(\frac{b_k(p)}{n}, X_{\fraca{b_k(p)}{n}}\biggr)  \sigma
_s^2 \,\mathrm{d}s \nonumber
\\&&\hphantom{\lim_{p \rightarrow\infty} \limsup_{n \rightarrow\infty} n^{\fraca
14}  \Biggl\{  \Biggl(\sum_{k=0}^{j_n(p)}  \biggl(}{}   + \int_{\fraca{c_k(p)}{n}}^{\fraca
{b_{k+1}(p)}{n}} \sigma_i^2\biggl(\frac{c_k(p)}{n}, X_{\fraca{c_k(p)}{n}}\biggr)  \sigma_s^2 \,\mathrm{d}s  \biggr)\\
&&\hphantom{\lim_{p \rightarrow\infty} \limsup_{n \rightarrow\infty} n^{\fraca
14}  \Biggl\{  \Biggl(}{} +
\int_{\frac{i_n(p)}{n}}^{1} \sigma_i^2\biggl(\frac
{i_n(p)}{n}, X_{\fraca{i_n(p)}{n}}\biggr)  \sigma_s^2 \,\mathrm{d}s  \Biggr)- C_i  \Biggr\} = 0.
\nonumber
\end{eqnarray}
As in the proof of Theorem~\ref{klord}, we have
%
\begin{eqnarray} \label{Taylor}
&&\sigma_i^2(s, X_s) - \sigma_i^2\biggl(\frac{b_k(p)}{n}, X_{\fraca{b_k(p)}{n}}\biggr)\nonumber
\\[-8pt]
\\[-8pt]
&& \quad = \frac{\partial}{\partial y} \sigma_i^2\bigl(s, X_{\fraca{b_k(p)}{n}}\bigr)
\biggl(\int_{\fraca{b_k(p)}{n}}^{s} \sigma_u \,\mathrm{d}W_u  \biggr)  + \mathrm{O}_p\biggl(\frac{p
m_n}{n} \biggr),
\nonumber
\end{eqnarray}
thus
%
\begin{eqnarray}\label{impl}
&&\int_{\fraca{b_k(p)}{n}}^{\fraca{c_k(p)}{n}}  \biggl(\sigma_i^2(s, X_s) -
\sigma_i^2\biggl(\frac{b_k(p)}{n}, X_{\fraca{b_k(p)}{n}}\biggr) \biggr)  \sigma_s^2 \,\mathrm{d}s
\nonumber
\\[-8pt]
\\[-8pt]&& \quad = \delta'(k,p)_3^n + \delta''(k,p)_3^n + \mathrm{O}_p\biggl(\frac{p^2
m_n^2}{n^2} \biggr),
\nonumber
\end{eqnarray}
where
\begin{eqnarray*}
\delta'(k,p)_3^n = \sigma_{\fraca{b_k(p)}{n}}^3 \int_{\fraca
{b_k(p)}{n}}^{\fraca{c_k(p)}{n}}
\frac{\partial}{\partial y} \sigma_i^2\bigl(s, X_{\fraca{b_k(p)}{n}}\bigr)
\biggl(\int_{\fraca{b_k(p)}{n}}^{s}  \mathrm{d}W_u  \biggr) \,\mathrm{d}s
\end{eqnarray*}
and $\delta''(k,p)_3^n$ is defined implicitly by equation (\ref{impl}).
We obtain
\begin{eqnarray*} \lim_{p \rightarrow\infty} \limsup_{n \rightarrow
\infty} n^{\fraca12} E \Biggl[ \Biggl (\sum_{k=0}^{j_n(p)} \delta
'(k,p)_3^n \Biggr)^2  \Biggr] = 0
\end{eqnarray*}
from the usual martingale argument and also
\begin{eqnarray*} \lim_{p \rightarrow\infty} \limsup_{n \rightarrow
\infty} n^{\fraca14} \sum_{k=0}^{j_n(p)} E[|\delta''(k,p)_3^n|]
\leq\lim_{p \rightarrow\infty} \limsup_{n \rightarrow\infty} K  p^{\fraca32}  n^{-\fraca12} \sum_{k=0}^{j_n(p)} E[\chi(p)_k^n] =
0
\end{eqnarray*}
as before. The corresponding results for the other summands in (\ref
{prop1}) can be shown analogously.

To finish the proof of Lemma~\ref{lem1}, we have to show
\begin{eqnarray*}
&&\lim_{p \rightarrow\infty} \limsup_{n \rightarrow\infty} n^{\fraca
14}  \Biggl\{ \sum_{k=0}^{j_n(p)}  \Biggl( \sigma_i^2\biggl(\frac{b_k(p)}{n},
X_{\fraca{b_k(p)}{n}}\biggr) \Biggl (\frac{1}{\kappa\psi_2} n^{-\fraca12} \sum
_{j=b_k(p)}^{c_k(p)-1} D(j)_3^n - \int_{\fraca{b_k(p)}{n}}^{\fraca
{c_k(p)}{n}} \sigma_s^2 \,\mathrm{d}s \Biggr) \\
 &&\hphantom{\lim_{p \rightarrow\infty} \limsup_{n \rightarrow\infty} n^{\fraca
14}  \Biggl\{ \sum_{k=0}^{j_n(p)}  \Biggl(}{}+ \sigma_i^2\biggl(\frac
{c_k(p)}{n}, X_{\fraca{c_k(p)}{n}}\biggr) \\
&&\hphantom{\lim_{p \rightarrow\infty} \limsup_{n \rightarrow\infty} n^{\fraca
14}  \Biggl\{ \sum_{k=0}^{j_n(p)}  \Biggl({}+}{}\times \Biggl(\frac{1}{\kappa\psi_2}
n^{-\fraca12} \sum_{j=c_k(p)}^{b_{k+1}(p)-1} D(j)_3^n - \int_{\fraca
{c_k(p)}{n}}^{\fraca{b_{k+1}(p)}{n}} \sigma_s^2 \,\mathrm{d}s  \Biggr)  \Biggr) \\
&&\hphantom{\lim_{p \rightarrow\infty} \limsup_{n \rightarrow\infty} n^{\fraca
14}  \Biggl\{ \sum_{k=0}^{j_n(p)}  }{}+
\sigma_i^2\biggl(\frac{i_n(p)}{n}, X_{\fraca{i_n(p)}{n}}\biggr)  \\
&&\hphantom{\lim_{p \rightarrow\infty} \limsup_{n \rightarrow\infty} n^{\fraca
14}  \Biggl\{ \sum_{k=0}^{j_n(p)}  {}+}{}\times\Biggl(\frac{1}{\kappa
\psi_2} n^{-\fraca12} \sum_{j=i_n(p)}^{n-m_n} D(j)_3^n - \int_{\fraca
{i_n(p)}{n}}^{1} \sigma_s^2 \,\mathrm{d}s  \Biggr)  \Biggr\} = 0,
\end{eqnarray*}
The last term is negligible, and the main idea for the tedious proof of
the remaining terms is to fix $k$ for a moment and to prove a
representation of the form
%
\begin{eqnarray} \label{interval}
\frac{1}{\kappa\psi_2} n^{-\fraca12} \sum_{j=b_k(p)}^{c_k(p)-1}
D(j)_3^n = \int_{\fraca{b_k(p)}{n}}^{\fraca{b_{k+1}(p)}{n}} h_{n,p}
\biggl(s-\frac{b_k(p)}{n} \biggr)  \sigma_s^2 \,\mathrm{d}s
\end{eqnarray}
for a suitable function $h_{n,p}(s)$, using the definition of
$D(j)_3^n$. A similar expression can be found for the sum from $c_k(p)$
to $b_{k+1}(p)$ with some $\bar h_{n,p} (s)$. A careful computation
shows that $h_{n,p}(s)$ is either close to one (for $s$ in the center
of the corresponding interval) or that $h_{n,p}(s)$ and $\bar
h_{n,p}(s)$ sum up to one (on its boundary). Then a Taylor expansion as
in the proof of \eqref{prop1} gives the result.
\end{pf}

\begin{lemma} \label{lem2}
We have
\begin{eqnarray*} \lim_{p \rightarrow\infty} \limsup_{n \rightarrow
\infty} n^{\fraca14} \Biggl \{ F_{in} -  \Biggl(\sum_{k=0}^{j_n(p)}
\bigl(G(k,p)_1^n +
G(k,p)_2^n\bigr) + G(p)_3^n  \Biggr)  \Biggr\} = 0.
\end{eqnarray*}
\end{lemma}
\begin{pf}
Without loss of generality, is suffices to show
\begin{eqnarray*} &&\lim_{p \rightarrow\infty} \limsup_{n \rightarrow
\infty} n^{-\fraca14} \sum_{k=0}^{j_n(p)} \sum_{j=b_k(p)}^{c_k(p)-1}
\biggl (\sigma_i^2(s, X_s) - \sigma_i^2\biggl(\frac{b_k(p)}{n},
X_{\fraca
{b_k(p)}{n}}\biggr) \biggr)  \biggl (|\overline{Z}^n_j|^2 - n^{-\fraca12} \frac
{\psi_1}{\kappa} \omega^2 \biggr)\\
&& \quad  = 0.
\end{eqnarray*}
The proof of this claim is tedious again. Essentially one simplifies
the expression above by the Taylor expansion from (\ref{Taylor}) and a
similar decomposition as in (\ref{ito}) for $|\overline{Z}^n_j|^2$ and
discusses each term separately.
\end{pf}

Note that we have completely analogous results for a decomposition of
$\hat{B}_t^{0} - B_t^0$. Thus, we end up with
%
\begin{eqnarray}
\label{mgda} \lim_{p
\rightarrow\infty} \limsup_{n \rightarrow\infty} n^{\fraca14} \Biggl \{
(\hat{B}_t^{0} - B_t^0) - \sum_{k=0}^{j_n(p)} Y(k,p) 1_{
\{\fraca{c_k(p)}{n} \leq t\}}  \Biggr\} &=& 0,
\\ \lim_{p \rightarrow\infty
} \limsup_{n \rightarrow\infty} n^{\fraca14}  \Biggl\{ (\hat C_i - C_i)
- \sum_{k=0}^{j_n(p)} \sigma_i^2\biggl(\frac{b_k(p)}{n}, X_{\fraca
{b_k(p)}{n}}\biggr)  Y(k,p)  \Biggr\} &=& 0  ,
\nonumber
\end{eqnarray}
where $Y(k,p)$ was defined in (\ref{defY}). Since
\begin{eqnarray*}
n  E\bigl[(Y(k,p))^2|\mathcal F_{\fraca{b_k(p)}{n}}\bigr] = p  \kappa \gamma
^2_{\fraca{b_k(p)}{n}} + \mathrm{o}_p(1) \quad\mbox{and} \quad
E\bigl[Y(k,p)|\mathcal F_{\fraca{b_k(p)}{n}}\bigr] = 0
\end{eqnarray*}
as in Jacod \textit{et al.}~\cite{jlmpv}, we conclude
\begin{eqnarray*}
&&\lim_{p \rightarrow\infty} \lim_{n \rightarrow\infty} n^{\fraca12}
\sum_{k=0}^{j_n(p)} E\bigl[Y(k,p)^2 1_{ \{\fraca{c_k(p)}{n} \leq t_i \wedge
t_j\}}|\mathcal F_{\fraca{b_k(p)}{n}}\bigr]\\
&& \quad   = \int_0^1 \gamma_s^2
 1_{[0,t_i \wedge t_j]}(s) \,\mathrm{d}s, \\
 &&\lim_{p
\rightarrow\infty} \lim_{n \rightarrow\infty} n^{\fraca12} \sum
_{k=0}^{j_n(p)} E\biggl[ Y(k,p)^2 1_{ \{\fraca{c_k(p)}{n} \leq t_i\}}
\sigma_i^2\biggl(\frac{b_k(p)}{n}, X_{\fraca{b_k(p)}{n}}\biggr) \big|\mathcal
F_{\fraca
{b_k(p)}{n}}\biggr]\\
&& \quad  = \int_0^1 \gamma_s^2  1_{[0,t_i]}(s)  \sigma
_j^2 (s, X_s) \,\mathrm{d}s, \\
&&\lim_{p \rightarrow\infty} \lim_{n \rightarrow\infty} n^{\fraca12}
\sum_{k=0}^{j_n(p)} E\biggl[Y(k,p)^2 \sigma_i^2\biggl(\frac{b_k(p)}{n},
X_{\fraca
{b_k(p)}{n}}\biggr) \sigma_j^2\biggl(\frac{b_k(p)}{n}, X_{\fraca{b_k(p)}{n}}\biggr)
\big|\mathcal F_{\fraca{b_k(p)}{n}}\biggr]
\\&& \quad = \int_0^1 \gamma_s^2  \sigma_i^2 (s, X_s)  \sigma_j^2 (s, X_s)
\,\mathrm{d}s.
\end{eqnarray*}
Theorem~\ref{fin} follows now from Theorem IX 7.28 in Jacod and
Shiryaev~\cite{jacshir}, since the missing
conditions can be shown in the same way as in Jacod \textit{et al.}
\cite{jlmpv}.\hfill\qed\vspace*{6pt}

The convergence of the finite dimensional distributions follows from
the delta method for stably converging sequences, since we have
\begin{eqnarray*}
n^{\fraca14} (\hat N_{t_1} - N_{t_1}, \ldots, \hat N_{t_k} - N_{t_k})^T
\stab Y \int_0^1 \Sigma^{\fraca12}_{t_1, \ldots, t_k} (s, X_s) \,\mathrm{d}W_s,
\end{eqnarray*}
where the $k \times(d+k)$-dimensional matrix $Y$ has the form
\begin{eqnarray*}
Y =
( I_{k \times k} \enskip -Y^*
)
, \qquad Y^* =
( B_{t_1}^T D^{-1} \enskip\cdots \enskip B_{t_k}^T D^{-1}
)
^T.
\end{eqnarray*}
A straightforward calculation shows that the conditional covariance
coincides with the one of the finite dimensional distributions of the
process defined in (\ref{dhat}). We are left to prove the tightness of
the process $n^{\fraca14} (\hat N_t - N_t)$, and this can be done by an
application of Theorem~VI.4.5 in Jacod and Shiryaev~\cite{jacshir},
using the boundedness of the processes involved as well as \mbox{$E[|\det
(D)|^{-\beta}] < \infty$}.
\end{pf}
\end{appendix}

\section*{Acknowledgements}
The authors would like to thank Martina Stein, who typed parts of this
manuscript with considerable technical expertise.
This work has been supported in part by the Collaborative
Research Center ``Statistical modeling of nonlinear dynamic processes''
(SFB~823) of the German Research Foundation (DFG).
The authors are also grateful
to the referees and the associate editor for constructive
comments on an earlier version of this paper.\looseness=1

%

\printhistory


\begin{thebibliography}{27}

\bibitem{ach}
%
\begin{bbook}[auto:STB|2011/11/17|08:29:20]
\bauthor{\bsnm{Achieser},~\bfnm{N.~J.}\binits{N.J.}}
(\byear{1956}).
\btitle{Theory of Approximation}.
\baddress{New York}: \bpublisher{Dover Publications Inc}.
\bptok{imsref}%
\end{bbook}
%
\endbibitem

\bibitem{ag}
%
\begin{barticle}[auto:STB|2011/11/17|08:29:20]
\bauthor{\bsnm{Ahn},~\bfnm{D.}\binits{D.}} \AND
\bauthor{\bsnm{Gao},~\bfnm{B.}\binits{B.}}
(\byear{1999}).
\btitle{A parametric nonlinear model of term structure dynamics}.
\bjournal{Review of Financial Studies}
\bvolume{12}
\bpages{721--762}.
\bptok{imsref}%
\end{barticle}
%
\endbibitem

\bibitem{ait}
%
\begin{barticle}[auto:STB|2011/11/17|08:29:20]
\bauthor{\bsnm{Ait-Sahalia},~\bfnm{Y.}\binits{Y.}}
(\byear{1996}).
\btitle{Testing continuous-time models of the spot interest rate}.
\bjournal{Review of Financial Studies}
\bvolume{9}
\bpages{385--426}.
\bptok{imsref}%
\end{barticle}
%
\endbibitem

\bibitem{am1987}
%
\begin{barticle}[auto:STB|2011/11/17|08:29:20]
\bauthor{\bsnm{Amihud},~\bfnm{Y.}\binits{Y.}} \AND
\bauthor{\bsnm{Mendelson},~\bfnm{Haim}\binits{H.}}
(\byear{1987}).
\btitle{Trading mechanisms and stock returns: An empirical investigation}.
\bjournal{J. Finance}
\bvolume{42}
\bpages{533--553}.
\bptok{imsref}%
\end{barticle}
%
\endbibitem

\bibitem{black1986}
%
\begin{barticle}[auto:STB|2011/11/17|08:29:20]
\bauthor{\bsnm{Black},~\bfnm{F.}\binits{F.}}
(\byear{1986}).
\btitle{Noise}.
\bjournal{J. Finance}
\bvolume{41}
\bpages{529--543}.
\bptok{imsref}%
\end{barticle}
%
\endbibitem

\bibitem{bs}
%
\begin{barticle}[auto:STB|2011/11/17|08:29:20]
\bauthor{\bsnm{Black},~\bfnm{F.}\binits{F.}} \AND
\bauthor{\bsnm{Scholes},~\bfnm{M.}\binits{M.}}
(\byear{1973}).
\btitle{The Pricing of Options and Corporate Liabilities}.
\bjournal{Journal of Political Economy}
\bvolume{81}
\bpages{637--659}.
\bptok{imsref}%
\end{barticle}
%
\endbibitem

\bibitem{chan}
%
\begin{barticle}[auto:STB|2011/11/17|08:29:20]
\bauthor{\bsnm{Chan},~\bfnm{K.~C.}\binits{K.C.}},
\bauthor{\bsnm{Karolyi},~\bfnm{G.~A.}\binits{G.A.}},
\bauthor{\bsnm{Longstaff},~\bfnm{F.~A.}\binits{F.A.}} \AND
\bauthor{\bsnm{Sanders},~\bfnm{A.~B.}\binits{A.B.}}
(\byear{1992}).
\btitle{An empirical~comparison of alternative models of the short-term
interest rate}.
\bjournal{J. Finance}
\bvolume{47}
\bpages{1209--1227}.
\bptok{imsref}%
\end{barticle}
%
\endbibitem

\bibitem{corr}
%
\begin{barticle}[mr]
\bauthor{\bsnm{Corradi},~\bfnm{Valentina}\binits{V.}} \AND
\bauthor{\bsnm{White},~\bfnm{Halbert}\binits{H.}}
(\byear{1999}).
\btitle{Specification tests for the variance of a diffusion}.
\bjournal{J.~Time Ser. Anal.}
\bvolume{20}
\bpages{253--270}.
\bid{doi={10.1111/1467-9892.00136}, issn={0143-9782}, mr={1693173}}
\bptok{imsref}%
\end{barticle}
%
\endbibitem

\bibitem{cir}
%
\begin{barticle}[mr]
\bauthor{\bsnm{Cox},~\bfnm{John~C.}\binits{J.C.}},
\bauthor{\bsnm{Ingersoll},~\bfnm{Jonathan~E.}\binits{J.E.} \bsuffix{Jr.}}
\AND\bauthor{\bsnm{Ross},~\bfnm{Stephen~A.}\binits{S.A.}}
(\byear{1985}).
\btitle{A theory of the term structure of interest rates}.
\bjournal{Econometrica}
\bvolume{53}
\bpages{385--407}.
\bid{doi={10.2307/1911242}, issn={0012-9682}, mr={0785475}}
\bptok{imsref}%
\end{barticle}
%
\endbibitem

\bibitem{dp}
%
\begin{barticle}[mr]
\bauthor{\bsnm{Dette},~\bfnm{Holger}\binits{H.}} \AND
\bauthor{\bsnm{Podolskij},~\bfnm{Mark}\binits{M.}}
(\byear{2008}).
\btitle{Testing the parametric form of the volatility in continuous time
diffusion models---a stochastic process approach}.
\bjournal{J. Econometrics}
\bvolume{143}
\bpages{56--73}.
\bid{doi={10.1016/j.jeconom.2007.08.002}, issn={0304-4076}, mr={2384433}}
\bptok{imsref}%
\end{barticle}
%
\endbibitem

\bibitem{dpv}
%
\begin{barticle}[mr]
\bauthor{\bsnm{Dette},~\bfnm{Holger}\binits{H.}},
\bauthor{\bsnm{Podolskij},~\bfnm{Mark}\binits{M.}} \AND
\bauthor{\bsnm{Vetter},~\bfnm{Mathias}\binits{M.}}
(\byear{2006}).
\btitle{Estimation of integrated volatility in continuous-time
financial models
with applications to goodness-of-fit testing}.
\bjournal{Scand. J. Statist.}
\bvolume{33}
\bpages{259--278}.
\bid{doi={10.1111/j.1467-9469.2006.00479.x}, issn={0303-6898}, mr={2279642}}
\bptok{imsref}%
\end{barticle}
%
\endbibitem

\bibitem{gallant1987}
%
\begin{bbook}[mr]
\bauthor{\bsnm{Gallant},~\bfnm{A.~Ronald}\binits{A.R.}}
(\byear{1987}).
\btitle{Nonlinear Statistical Models}.
\bseries{Wiley Series in Probability and Mathematical Statistics: Applied
Probability and Statistics}.
\baddress{New York}: \bpublisher{Wiley}.
\bid{doi={10.1002/9780470316719}, mr={0921029}}
\bptok{imsref}%
\end{bbook}
%
\endbibitem

\bibitem{glot}
%
\begin{barticle}[mr]
\bauthor{\bsnm{Gloter},~\bfnm{Arnaud}\binits{A.}} \AND
\bauthor{\bsnm{Jacod},~\bfnm{Jean}\binits{J.}}
(\byear{2001}).
\btitle{Diffusions with measurement errors. {II}. {O}ptimal estimators}.
\bjournal{ESAIM Probab. Statist.}
\bvolume{5}
\bpages{243--260 (electronic)}.
\bid{doi={10.1051/ps:2001111}, issn={1292-8100}, mr={1875673}}
\bptok{imsref}%
\end{barticle}
%
\endbibitem

\bibitem{harris1990}
%
\begin{barticle}[auto:STB|2011/11/17|08:29:20]
\bauthor{\bsnm{Harris},~\bfnm{L.}\binits{L.}}
(\byear{1990}).
\btitle{Estimation of stock variance and serial covariance from discrete
observations}.
\bjournal{Journal of Financial and Quantitative Analysis}
\bvolume{25}
\bpages{291--306}.
\bptok{imsref}%
\end{barticle}
%
\endbibitem

\bibitem{harris1991}
%
\begin{barticle}[auto:STB|2011/11/17|08:29:20]
\bauthor{\bsnm{Harris},~\bfnm{L.}\binits{L.}}
(\byear{1991}).
\btitle{Stock price clustering and discreteness}.
\bjournal{Review of Financial Studies}
\bvolume{4}
\bpages{389--415}.
\bptok{imsref}%
\end{barticle}
%
\endbibitem

\bibitem{heston}
%
\begin{barticle}[auto:STB|2011/11/17|08:29:20]
\bauthor{\bsnm{Heston},~\bfnm{S.~L.}\binits{S.L.}}
(\byear{1993}).
\btitle{A Closed-Form Solution for Options with Stochastic Volatility with
Applications to Bond and Currency Options}.
\bjournal{Review of Financial Studies}
\bvolume{6}
\bpages{327--343}.
\bptok{imsref}%
\end{barticle}
%
\endbibitem

\bibitem{hull}
%
\begin{barticle}[auto:STB|2011/11/17|08:29:20]
\bauthor{\bsnm{Hull},~\bfnm{J.}\binits{J.}} \AND
\bauthor{\bsnm{White},~\bfnm{A.}\binits{A.}}
(\byear{1987}).
\btitle{The Pricing of Options on Assets with Stochastic Volatilities}.
\bjournal{J. Finance}
\bvolume{42}
\bpages{281--300}.
\bptok{imsref}%
\end{barticle}
%
\endbibitem

\bibitem{jac08}
%
\begin{barticle}[mr]
\bauthor{\bsnm{Jacod},~\bfnm{Jean}\binits{J.}}
(\byear{2008}).
\btitle{Asymptotic properties of realized power variations and related
functionals of semimartingales}.
\bjournal{Stochastic Process. Appl.}
\bvolume{118}
\bpages{517--559}.
\bid{doi={10.1016/j.spa.2007.05.005}, issn={0304-4149}, mr={2394762}}
\bptok{imsref}%
\end{barticle}
%
\endbibitem

\bibitem{jlmpv}
%
\begin{barticle}[mr]
\bauthor{\bsnm{Jacod},~\bfnm{Jean}\binits{J.}},
\bauthor{\bsnm{Li},~\bfnm{Yingying}\binits{Y.}},
\bauthor{\bsnm{Mykland},~\bfnm{Per~A.}\binits{P.A.}},
\bauthor{\bsnm{Podolskij},~\bfnm{Mark}\binits{M.}} \AND
\bauthor{\bsnm{Vetter},~\bfnm{Mathias}\binits{M.}}
(\byear{2009}).
\btitle{Microstructure noise in the continuous case: The pre-averaging
approach}.
\bjournal{Stochastic Process. Appl.}
\bvolume{119}
\bpages{2249--2276}.
\bid{doi={10.1016/j.spa.2008.11.004}, issn={0304-4149}, mr={2531091}}%
\bptnote{check year}%
\bptok{imsref}%
\end{barticle}
%
\endbibitem

\bibitem{jacshir}
%
\begin{bbook}[mr]
\bauthor{\bsnm{Jacod},~\bfnm{Jean}\binits{J.}} \AND
\bauthor{\bsnm{Shiryaev},~\bfnm{Albert~N.}\binits{A.N.}}
(\byear{2003}).
\btitle{Limit Theorems for Stochastic Processes},
\bedition{2nd} ed.
\bseries{Grundlehren der Mathematischen Wissenschaften [Fundamental Principles
of Mathematical Sciences]}
\bvolume{288}.
\baddress{Berlin}: \bpublisher{Springer}.
\bid{mr={1943877}}
\bptok{imsref}%
\end{bbook}
%
\endbibitem

\bibitem{pv2}
%
\begin{barticle}[mr]
\bauthor{\bsnm{Podolskij},~\bfnm{Mark}\binits{M.}} \AND
\bauthor{\bsnm{Vetter},~\bfnm{Mathias}\binits{M.}}
(\byear{2009}).
\btitle{Bipower-type estimation in a noisy diffusion setting}.
\bjournal{Stochastic Process. Appl.}
\bvolume{119}
\bpages{2803--2831}.
\bid{doi={10.1016/j.spa.2009.02.006}, issn={0304-4149}, mr={2554029}}
\bptok{imsref}%
\end{barticle}
%
\endbibitem

\bibitem{pv1}
%
\begin{barticle}[mr]
\bauthor{\bsnm{Podolskij},~\bfnm{Mark}\binits{M.}} \AND
\bauthor{\bsnm{Vetter},~\bfnm{Mathias}\binits{M.}}
(\byear{2009}).
\btitle{Estimation of volatility functionals in the simultaneous
presence of
microstructure noise and jumps}.
\bjournal{Bernoulli}
\bvolume{15}
\bpages{634--658}.
\bid{doi={10.3150/08-BEJ167}, issn={1350-7265}, mr={2555193}}
\bptok{imsref}%
\end{barticle}
%
\endbibitem

\bibitem{rev}
%
\begin{bbook}[mr]
\bauthor{\bsnm{Revuz},~\bfnm{Daniel}\binits{D.}} \AND
\bauthor{\bsnm{Yor},~\bfnm{Marc}\binits{M.}}
(\byear{1999}).
\btitle{Continuous Martingales and {B}rownian Motion},
\bedition{3rd} ed.
\bseries{Grundlehren der Mathematischen Wissenschaften [Fundamental Principles
of Mathematical Sciences]}
\bvolume{293}.
\baddress{Berlin}: \bpublisher{Springer}.
\bid{mr={1725357}}
\bptok{imsref}%
\end{bbook}
%
\endbibitem

\bibitem{sw1989}
%
\begin{bbook}[mr]
\bauthor{\bsnm{Seber},~\bfnm{G.~A.~F.}\binits{G.A.F.}} \AND
\bauthor{\bsnm{Wild},~\bfnm{C.~J.}\binits{C.J.}}
(\byear{1989}).
\btitle{Nonlinear Regression}.
\bseries{Wiley Series in Probability and Mathematical Statistics: Probability
and Mathematical Statistics}.
\baddress{New York}: \bpublisher{Wiley}.
\bid{doi={10.1002/0471725315}, mr={0986070}}
\bptok{imsref}%
\end{bbook}
%
\endbibitem

\bibitem{vas}
%
\begin{barticle}[auto:STB|2011/11/17|08:29:20]
\bauthor{\bsnm{Vasicek},~\bfnm{O.}\binits{O.}}
(\byear{1977}).
\btitle{An equilibrium characterization of the term structure}.
\bjournal{Journal of Financial Economics}
\bvolume{5}
\bpages{177--188}.
\bptok{imsref}%
\end{barticle}
%
\endbibitem

\bibitem{zha}
%
\begin{barticle}[mr]
\bauthor{\bsnm{Zhang},~\bfnm{Lan}\binits{L.}},
\bauthor{\bsnm{Mykland},~\bfnm{Per~A.}\binits{P.A.}} \AND
\bauthor{\bsnm{A{\"{\i}}t-Sahalia},~\bfnm{Yacine}\binits{Y.}}
(\byear{2005}).
\btitle{A tale of two time scales: Determining integrated volatility
with noisy
high-frequency data}.
\bjournal{J. Amer. Statist. Assoc.}
\bvolume{100}
\bpages{1394--1411}.
\bid{doi={10.1198/016214505000000169}, issn={0162-1459}, mr={2236450}}
\bptok{imsref}%
\end{barticle}
%
\endbibitem

\end{thebibliography}
\end{document}